\documentclass[11pt]{amsart}

\usepackage{amsmath,amssymb,amsthm,mathtools}
\DeclareMathOperator*{\argmin}{arg\,min}
\usepackage{bm}

\usepackage[T1]{fontenc}
\usepackage[utf8]{inputenc}
\usepackage{float}
\usepackage{tikz}
\usetikzlibrary{matrix,arrows.meta,positioning,calc,decorations.pathreplacing,shapes.geometric,shapes.misc,fit,backgrounds,patterns}
\usepackage{theoremref}
\usepackage{lmodern}
\usepackage{microtype}
\usepackage[a4paper,margin=3cm]{geometry}

\usepackage{enumitem}
\usepackage{booktabs}
\usepackage{graphicx}

\usepackage{xcolor}
\usepackage[
    colorlinks=true,
    linkcolor=blue,
    citecolor=blue,
    urlcolor=blue
]{hyperref}

\newtheorem{theorem}{Theorem}[section]
\newtheorem{proposition}[theorem]{Proposition}
\newtheorem{lemma}[theorem]{Lemma}

\theoremstyle{definition}
\newtheorem{definition}[theorem]{Definition}
\newtheorem{example}[theorem]{Example}

\theoremstyle{remark}
\newtheorem{remark}[theorem]{Remark}

\newcommand{\R}{\mathbb{R}}

\newcommand{\Z}{\mathbb{Z}}
\newcommand{\N}{\mathbb{N}}
\newcommand{\cA}{\mathcal{A}}
\newcommand{\cS}{\mathcal{S}}
\newcommand{\cE}{\mathcal{E}}
\newcommand{\cQ}{\mathcal{Q}}
\newcommand{\cR}{\mathcal{R}}
\DeclareMathOperator{\Pos}{Pos}

\newcommand{\supp}{\operatorname{supp}}
\newcommand{\conv}{\operatorname{conv}}

\title[TMP on $\Z$]{The truncated $\Z-$moment problem}

\author{Maria Infusino}
\address{
Department of Mathematics and Computer Science\\
University of Cagliari\\
Italy
}
\email{maria.infusino@unica.it}

\author{Luca Onnis}
\address{
Department of Mathematics\\
University of Genova\\
Italy
}
\email{luca.onnis@edu.unige.it}

\date{\today}



\begin{document}

\begin{abstract}
We study the truncated moment problem on the one-dimensional full integer
lattice $\Z$, seeking necessary and sufficient conditions
for a prescribed sequence of $n$ real numbers to be $\Z-$representable, i.e. to be the sequence of the first
$n$ moments of a probability measure supported on $\Z$.
Building on the approach developed in
\cite{InfusinoKunaLebowitzSpeer2017} for the truncated moment problem on the non-negative integers, we establish a
complete characterization of $\Z-$representable sequences for every truncation order $n$. From this result, we derive an inductive and computable $\Z-$representability criterion. More generally, we also investigate the relationship between the truncated $\Z-$moment problem and the truncated Hamburger moment problem, seen as its continuous counterpart.
\end{abstract}
\maketitle
\vspace{-5mm}
\section{Introduction}
The truncated moment problem asks whether finitely many prescribed real numbers
can be represented as the moments of a probability measure with a specified
support. When the support is fixed to be the real line, the problem is classically called the Hamburger moment problem and the solvability of its truncated version is characterized by positivity and
rank conditions on the associated Hankel matrix~\cite{Hamburger1920,ShohatTamarkin1943,CurtoFialkow1991}
and by the existence of flat extensions~\cite{CurtoFialkow1996,Schmudgen2017}.
Requiring integer support imposes additional constraints that are not captured by the positivity of the Hankel matrix. For example, if the mean is $1/2$, the smallest possible second moment is $1/4$ on $\mathbb{R}$, but $1/2$ on $\mathbb{Z}$, since $x(x-1)\ge 0$ for every $x\in\mathbb{Z}$ forces the second moment to be greater than or equal to the first.

Integer-valued random variables naturally arise in many applied settings, such as the
comparison of counts from different populations~\cite{Irwin1937,Skellam1946}, the
modelling of score differences in sports~\cite{Maher1982,KarlisNtzoufras2009},
integer-valued time series~\cite{AlOshAlzaid1987,Freeland2010}, and the analysis of
tick-by-tick price changes in finance~\cite{BarndorffNielsenPollardShephard2012}. For instance, the parameters of the Skellam distribution~\cite{Skellam1946}, which
models the difference of two independent Poisson counts, are determined by its mean
and variance alone, so checking whether prescribed low-order moments are compatible
with an integer-valued distribution is precisely a truncated $\Z-$moment problem. Moment problems on discrete sets also play an important role in optimization,
for instance in the linear programming approach to bounding probabilities from known
moments~\cite{Prekopa1990} and in the derivation of optimal probability inequalities
via convex optimization~\cite{BertsimasPopescu2005}.

The set of vectors representable by discretely supported probability measures has been studied using convex geometry and
linear programming (see e.g. \cite{KarlinStudden1966,Kemp1968, KN1977, Prekopa1990}). In particular, for the non-negative integers $\mathbb{N}_0$, Infusino, Kuna, Lebowitz
and Speer~\cite{InfusinoKunaLebowitzSpeer2017} obtained the first computable necessary
and sufficient conditions for a vector of $n$ real numbers to be the sequence of the
first $n$ moments of a probability measure supported on $\mathbb{N}_0$, for every
truncation order~$n$. Their approach is an inductive procedure that parallels the
truncated Stieltjes moment problem and distinguishes interior and boundary
representable moment vectors.  Here we consider the full integer
lattice $\mathbb{Z}$, which is unbounded in both directions and therefore falls
outside that framework. Our main results are a solvability criterion for the truncated $\Z-$moment problem, for every truncation order $n$, and a finite recursive algorithm for deciding it.

Let us describe in detail the structure of this article. In Section~\ref{sec:preliminaries}
we introduce the essential notation and recall some basic results on the truncated $K-$moment
problem. In Section~\ref{sec:finite-lattice-Z} we provide a solvability criterion for the moment problem on a finite integer interval
(see \thref{thm:Z-finite-support}). In Section~\ref{sec:Z-moment-problem} we extend this analysis to the full integer lattice~$\Z$ and we establish the main result: a solvability criterion for the $\Z-$moment problem
(see \thref{thm:Z-main}). In Section~\ref{sec:extremal-completions-Z} we develop an inductive
solvability criterion (see \thref{prop:Z-inductive}), we compare the
extremal $\Z-$ and Hamburger completions (see \thref{def:extremal-Z-completion} and \thref{def:Hamburger-completion}) and we prove an interlacing
property of their supports (see \thref{prop:Z-hamburger-interlacing}).
In Section~\ref{sec:decision-Z} we present a finite recursive algorithm
for deciding $\Z-$representability, practically illustrated in
\thref{ex:Z-algorithm-order-four}. Finally, the auxiliary technical results are collected in Appendix~\ref{app:core-technical-results}.

\section{Preliminaries}
\label{sec:preliminaries}

Throughout the paper, we deal with probability measures so we fix the zero-th moment $m_0=1$. Let us recall the statement of the $K-$moment problem and some related basic definitions.
\begin{definition}[$n-$truncated $K-$moment problem]
\label{def:truncated-K-moment}
Let $n\in\N$ and $K\subseteq\R$ be closed and non-empty.  Given $m=(m_1,\ldots,m_n)\in\R^n$, we ask if there is a probability measure $\mu$ with $\operatorname{supp}{\mu}\subseteq K$ such that
\[
 m_k=\int_K x^k\,d\mu(x),\qquad \forall k\in\{1,\ldots,n\}.
\]
Such a measure is called a \emph{$K-$representing probability measure for $m$}.
When such a measure exists, $m$ is called \emph{$K-$representable}, or an
$n-$truncated \emph{$K-$moment vector}. We further denote the set of all such vectors by $\mathcal{S}_n(K)$.
\end{definition}
For $m=(m_1,\ldots,m_n)\in\R^n$ and $t\in\R$, we set
$(m,t):=(m_1,\ldots,m_n,t)$.
We further associate with \(m=(m_1,\ldots,m_n)\in\R^n\) the corresponding Riesz functional defined by
\[
\begin{aligned}
L_m:\mathbb{R}[x]_{\le n} &\longrightarrow \mathbb{R},\\
\sum_{k=0}^{n} p_k x^k &\longmapsto p_0 + \sum_{k=1}^{n} p_k m_k,
\end{aligned}
\]
where $\R[x]_{\le n}$ denotes the polynomials of degree at most $n$.
The Riesz functional provides a necessary condition for the solvability of the \(n\)-truncated \(K\)-moment problem (cf. \thref{prop:positivity-support}-(i)). To formulate this condition, we first need to introduce some notation.

\begin{definition}[Cones of non-negative polynomials on $K$]
\label{def:positive-cones}
For $K\subseteq\R$, we define
\[
\begin{aligned}
 \Pos(K)&:=\{p\in\R[x]:p(t)\ge0\text{ }\forall t\in K\},
\end{aligned}
\]
and we denote its truncation to degree at most $n$ by $\Pos(K)_{\le n}:=\Pos(K)\cap\R[x]_{\le n}.$
\end{definition}

For $p\in\R[x]$, we write
\[
Z(p):=\{t\in\R:p(t)=0\},\qquad Z_K(p):=Z(p)\cap K.
\]
We will repeatedly use the following support restriction and atomic reduction properties.

\begin{proposition}[Positivity and support restriction]
\thlabel{prop:positivity-support}
Let $m\in\cS_n(K)$ and $\mu$ be a $K-$representing probability measure for $m$. Then
\begin{enumerate}[label=\textup{(\roman*)}]
\item $L_m(p)\ge0$ for every $p\in\Pos(K)_{\le n}$.
\item If $p\in\Pos(K)_{\le n}$ and $L_m(p)=0$, then $\supp\mu\subseteq Z_K(p)$.
\end{enumerate}
\end{proposition}

\begin{proof}
Since $L_m(p)=\int_Kp\,d\mu$, assertion \textup{(i)} follows from
non-negativity. If this integral is zero, then $p=0$ $\mu$-almost
everywhere and, therefore, every point of $\supp\mu$ is a zero of $p$, which proves \textup{(ii)}.
\end{proof}

By applying the finite-atomic reduction of Richter, Rogosinski and Tchakaloff (see \cite{Richter1957,Rogosinski1958, Tchakaloff1957}) to the functions $1,x,\ldots,x^n$ we obtain the following statement.

\begin{proposition}[Atomic reduction]
\thlabel{prop:atomic-reduction}
If $m\in\cS_n(K)$, then $m$ has a $K-$representing probability measure supported on at most $n+1$ points. In particular,
\[
\begin{aligned}
 m\in\cS_n(\Z)
 &\quad\Longleftrightarrow\quad
 m\in\cS_n(\{-N,\ldots,N\})\text{ for some }N\in\N.
\end{aligned}
\]

\end{proposition}

Convex geometry provides a natural framework for the study of the truncated $K-$moment problem. In particular, we show in the next proposition that if $K$ is finite then $\cS_n(K)$ is a polytope. To this aim, let us introduce the following concept.
\begin{definition}[Moment curves]
\thlabel{def:momentcurve}
Let \(n\in\mathbb{N}\) and $K\subseteq \R$. The \emph{moment curve} of order \(n\) is
\[
\gamma_n:K\to\mathbb{R}^{n+1},\qquad
\gamma_n(x):=(1,x,x^2,\dots,x^n).
\]
Its projection onto the last \(n\) coordinates given by
\[
v_n(x):=(x,x^2,\dots,x^n)=\pi_n(\gamma_n(x)),
\]
where \(\pi_n(x_0,\dots,x_n):=(x_1,\dots,x_n)\), is called the
\emph{projected moment curve} of order \(n\).
\end{definition}
\begin{proposition}
\thlabel{prop:finite-support-polytope}
Let $n\in\N$ and $M\in\N_0$. If $K=\{x_0,\dots,x_M\}\subset\R$, then
\[
\cS_n(K)=\conv\{v_n(x_0),\dots,v_n(x_M)\}.
\]
\end{proposition}

\begin{proof}
If $\mu$ is a probability measure on $\mathbb{R}$ supported on $K$, then it has the form $\mu=\sum_{j=0}^M \lambda_j\delta_{x_j},$ for some $\lambda_j\ge 0$ with $\sum_{j=0}^M\lambda_j=1$. Hence, its moment vector of order $n$ is $\sum_{j=0}^M\lambda_jv_n(x_j)$ and, thus, $\cS_n(K)\subseteq \conv\{v_n(x_0),\dots,v_n(x_M)\}.$
Conversely, let $y=\sum_{j=0}^M\lambda_jv_n(x_j),$ with $\lambda_j\ge 0$ for every $j\in\{0,\dots,M\}$ and $\sum_{j=0}^{M}\lambda_j=1.$
Then $y$ is represented by the probability measure $\mu=\sum_{j=0}^M\lambda_j\delta_{x_j}$
supported on \(K\). Hence,
\(\conv\{v_n(x_0),\dots,v_n(x_M)\}\subseteq \cS_n(K)\).
\end{proof}
\section{\texorpdfstring{The truncated $\{-N,\ldots,N\}-$moment problem}{The truncated \{-N,...,N\}-moment problem}}
\label{sec:finite-lattice-Z}
As a first step towards the solvability of the $n-$truncated $\Z-$moment problem, we analyze the $n-$truncated $\{-N,\ldots,N\}-$moment problem, where $N\in\N$, and then pass from these finite truncations to the full integer lattice.
By \thref{prop:finite-support-polytope}, this is a polyhedral membership problem: the
moment polytope $\cS_n(\{-N,\ldots,N\})$ is the convex hull of the images of the points in $\{-N,\ldots,N\}$ under the projected moment curve of order~$n$ introduced in \thref{def:momentcurve}. In view of this observation, we recall the following standard geometric lemma.
\begin{lemma}
\thlabel{lem:convex-hull-halfspaces}
Let $S\subset\R^n$ be finite with $\operatorname{aff}(S)=\R^n$.
Then
\begin{equation}
    \label{eq:conv-hull}
    \operatorname{conv}(S)=\bigcap_{H\in\mathcal H} H,
\end{equation}
where $\mathcal H$ is the family of all closed half-spaces $H$ containing $S$
whose bounding hyperplane $\partial H$ contains $n$ affinely independent
points of $S$.
Moreover, this representation is minimal: no half-space may be omitted from the
intersection in \eqref{eq:conv-hull}.
\end{lemma}

\begin{proof}
This is a consequence of \cite[Theorem~3.1.1]{Grunbaum2003ConvexPolytopes}.
\end{proof}
In order to derive solvability conditions for the $n-$truncated moment problem on
$K_N:=\{-N,\ldots,N\}$, we use \thref{lem:convex-hull-halfspaces}
to reduce the question to the characterization of the polynomials of degree $n$
that are non-negative on $K_N$ and have $n$ distinct zeros in this finite set. We describe these polynomials by their roots, using a notation that will
also allow us to treat the analogous class of polynomials non-negative
on the full integer lattice.
\begin{definition}
\thlabel{def:A_n}
We set \(\mathcal A_0:=\{\emptyset\}\) and \(p_{\emptyset}:=1\). For \(n\in\mathbb{N}\), we define \(\mathcal A_n\) as follows:
\begin{enumerate}[label=(\roman*),leftmargin=1.2cm]
\item if \(n\) is even, then \(\mathcal A_n\) is the set of \(n\)-tuples $\alpha=(\alpha_1,\dots,\alpha_n)\in\Z^n$ such that $\alpha_1<\alpha_2<\cdots<\alpha_n$ and $\alpha_{2k}=\alpha_{2k-1}+1$ for every $k\in\{1,\dots,n/2\};$
\item if \(n\) is odd, then $\mathcal A_n:=\varnothing.$
\end{enumerate}
For $\alpha=(\alpha_1,\dots,\alpha_n)\in\cA_n$ we define
\[
p_\alpha(x):=\prod_{i=1}^n(x-\alpha_i),
\]
and set $\cE_n:=\{p_\alpha:\alpha\in\mathcal A_n\}$. We also set \(\cE_0:=\{1\}\).
\end{definition}
Let us fix a useful short notation for certain polynomials often occurring in this work.
\begin{definition}[Lattice blocks]
\thlabel{def:block-poly}
    We define the \emph{lattice block} associated with $a\in\Z$ as $B_a(x) = (x-a)(x-a-1)$.
\end{definition}

   We can rewrite the
    family of polynomials \(\cE_n\) of \thref{def:A_n}
    in terms of lattice blocks.
    \begin{itemize}
        \item If \(n=2r\), then
        \begin{equation}
\label{eq:E-lattice-blocks}
\cE_{2r}
=
\left\{
\prod_{j=1}^{r}B_{a_j}(x):
a_j\in\Z,\;
a_{j+1}\ge a_j+2,\ \forall j\in\{1,\dots,r-1\}
\right\}.
\end{equation}
        \item If \(n=2r+1\), then $\cE_{2r+1}=\emptyset$.
    \end{itemize}

\begin{proposition}
\thlabel{prop:Z-structure}
For every \(n\in\mathbb{N}\), $\cE_n$ is the set of all monic polynomials of degree $n$ with $n$ distinct roots in $\Z$ which are non-negative on $\Z$.
\end{proposition}

\begin{proof}
Let \(p\) be a monic polynomial of degree \(n\) non-negative on $\Z$ with \(n\) distinct roots $r_1<r_2<\cdots<r_n$ in \(\Z\).
Because all roots are simple and \(p\) is monic, the sign of \(p\) alternates across consecutive roots. Furthermore, $p$ cannot have odd degree since if this were the case then $p$ would tend to $-\infty$ as $x\to-\infty$ and, thus, it would not be non-negative on~$\Z$. Therefore, if $n$ is odd, then $\mathcal{A}_n=\varnothing$ and $\cE_n=\varnothing$. From now on, we assume that $n$ is even. Any interval on which \(p\) is negative cannot contain an
integer, since by assumption \(p(j)\ge 0\) for all \(j\in\Z\).
Hence, if $p$ is negative on \((r_i,r_{i+1})\) for some $i\in\{1,\dots,n-1\}$, then this interval cannot contain any point of \(\Z\),
which forces $r_{i+1}=r_i+1.$ Since \(n\) is even, $p$ is non-negative on \((-\infty,r_1)\) and so all the intervals on which $p$ is negative are of the form $(r_{2k-1},r_{2k-1}+1)$ for some $k\in\{1,\dots,\frac{n}{2}\}$, giving \textup{(i)} in \thref{def:A_n}.
The converse inclusion is immediate, as for each $\alpha=(\alpha_1,\dots,\alpha_n)\in\cA_n$, every factor $(x-\alpha_{2k-1})(x-\alpha_{2k})=B_{\alpha_{2k-1}}(x)$ is non-negative on \(\Z\), and so is their product which gives~$p_{\alpha}$.
\end{proof}

\begin{definition}
\thlabel{def:four-families}
For \(N,n\in\mathbb N\) with \(2N\ge n\), we define the following four families of polynomials on $K_N$ with degree $n$:
\[
\cE_{n,N}:=
\{p_\alpha:\alpha\in\mathcal A_n,\ \alpha_1\ge -N,\ \alpha_n\le N\},
\]
\[
\cQ_{n,N}^{-}:=
\{(N+x)p_\alpha:\alpha\in\mathcal A_{n-1},\ \alpha_1\ge -N+1,\ \alpha_{n-1}\le N\},
\]
\[
\cQ_{n,N}^{+}:=
\{(N-x)p_\alpha:\alpha\in\mathcal A_{n-1},\ \alpha_1\ge -N,\ \alpha_{n-1}\le N-1\},
\]
and, for $n\ge2$,
\[
\cR_{n,N}:=
\{(N+x)(N-x)p_\alpha:\alpha\in\mathcal A_{n-2},\ \alpha_1\ge -N+1,\ \alpha_{n-2}\le N-1\},
\]
while $\cR_{1,N}:=\varnothing$. Here $p_\alpha$ and $\cA_n$ are as in \thref{def:A_n}, and whenever $\alpha=\emptyset$ the inequalities involving its coordinates are understood to be void.
\end{definition}

\begin{lemma}
\thlabel{lem:KN-zero-classification}
Let \(N,n\in\mathbb N\) with \(2N\ge n\), and let
\(q\in\Pos(K_N)_{\le n}\) be a polynomial of degree \(n\) with \(n\) distinct zeros
in \(K_N\). Then \(q\) is a positive scalar multiple of exactly one polynomial in one of the four families: $\cE_{n,N},
\cQ_{n,N}^{-},
\cQ_{n,N}^{+},
\cR_{n,N}.$
More explicitly:
\begin{enumerate}[label=(\roman*),leftmargin=1.2cm]
\item if \(q\) has positive leading coefficient and \(n\) is even, then
\(q=cp_\alpha\) for some \(c>0\) and some \(\alpha\in\mathcal A_n\) with
\(\alpha_1\ge -N\) and \(\alpha_n\le N\), i.e.\ $q/c\in\cE_{n,N}$.

\item if \(q\) has positive leading coefficient and \(n\) is odd, then
\(q=c(N+x)p_\alpha\) for some \(c>0\) and some
\(\alpha\in\mathcal A_{n-1}\) with \(\alpha_1\ge -N+1\) and
\(\alpha_{n-1}\le N\) when $n \geq 3$; if $n=1$, then $q(x)=c(N+x)$ for some $c>0$. In both cases $q/c\in\cQ_{n,N}^{-}$.

\item if \(q\) has negative leading coefficient and \(n\) is odd, then
\(q=c(N-x)p_\alpha\) for some \(c>0\) and some
\(\alpha\in\mathcal A_{n-1}\) with \(\alpha_1\ge -N\) and
\(\alpha_{n-1}\le N-1\) when $n \geq 3$; if $n=1$, then $q(x)=c(N-x)$ for some $c>0$. In both cases $q/c\in\cQ_{n,N}^{+}$.

\item if \(q\) has negative leading coefficient and \(n\) is even, then
\(q=c(N+x)(N-x)p_\alpha\) for some \(c>0\) and some
\(\alpha\in\mathcal A_{n-2}\) with \(\alpha_1\ge -N+1\) and
\(\alpha_{n-2}\le N-1\) when $n \geq 4$; if $n=2$, then $q(x)=c(N-x)(N+x)$ for some $c>0$. In both cases $q/c\in\cR_{n,N}$.
\end{enumerate}
\end{lemma}
\begin{proof}
For $n=1$, the conclusion follows directly from the fact that the unique zero must be one of the endpoints $-N$ or $N$, according to the sign of the leading coefficient. We may therefore assume that $n\geq2$. Suppose $q\in\Pos(K_N)_{\leq n}$ and write the zeros of \(q\) as
$-N\le r_1<r_2<\cdots<r_n\le N.$
Since \(\deg q=n\), we may write
\[
q(x)=c\prod_{i=1}^n(x-r_i),
\quad \mbox{with } c\in\mathbb{R}\setminus\{0\}.
\]
Assume first that \(c>0\). Then \(q(x)>0\) for all \(x>r_n\), and the sign of \(q\)
alternates when crossing a simple root. Since \(q\ge0\) on \(K_N\), every
interval on which \(q\) is negative must contain no point of \(K_N\) and, hence, any
two roots which bound such an interval have to be consecutive integers in \(K_N\).
If \(n\) is even, \(q\) is also positive on \((-\infty,r_1)\), and all the
intervals on which \(q\) is negative are intervals between two consecutive roots.
Therefore, $r_{2k}=r_{2k-1}+1$ for every $k\in\{1,\dots,n/2\}.$
Thus, \((r_1,\dots,r_n)\in\mathcal A_n\) with \(r_1\ge -N\) and
\(r_n\le N\) and, hence, $q=cp_\alpha$
for some \(c>0\) and some \(\alpha\in\mathcal A_n\) satisfying the required
bounds.
If \(n\) is odd, then \(q\) is negative on the interval
\((-\infty,r_1)\). Since this interval cannot contain any point of \(K_N\), we
must have \(r_1=-N\); otherwise \(r_1 > -N\) and
\(q(-N)<0\), which contradicts \(q\in\Pos(K_N)_{\leq n}\). Therefore, we get
\[
q(x)
=
c(N+x)\prod_{i=2}^n(x-r_i).
\]
The remaining intervals on which \(q\) is negative are bounded by consecutive
pairs of roots. Hence, $r_{2k+1}=r_{2k}+1$ for every $k\in\{1,\dots,(n-1)/2\}$,
and so \((r_2,\dots,r_n)\in\mathcal A_{n-1}\) with
\(r_2\ge -N+1\) and \(r_n\le N\). Consequently, $q=c(N+x)p_\alpha$
for some \(c>0\) and some \(\alpha\in\mathcal A_{n-1}\) satisfying the
required bounds.

Assume now that \(c<0\). Then \(q(x)<0\) for all \(x>r_n\) and, hence,
\(r_n\) must be the right endpoint \(N\); otherwise \(r_n < N\) and
\(q(N)<0\), contradicting \(q\in\Pos(K_N)_{\leq n}\).
If \(n\) is odd, \(q\) is positive on \((-\infty,r_1)\), and all the remaining
intervals on which \(q\) is negative are bounded by consecutive pairs of roots.
Therefore, $r_{2k}=r_{2k-1}+1$ for every $k=1,\dots,(n-1)/2.$
Moreover,
\[
q(x)
=
c(x-N)\prod_{i=1}^{n-1}(x-r_i)
=
(-c)(N-x)\prod_{i=1}^{n-1}(x-r_i).
\]
Thus, \((r_1,\dots,r_{n-1})\in\mathcal A_{n-1}\), with
\(r_1\ge -N\) and \(r_{n-1}\le N-1\) and, hence, $q=c'(N-x)p_\alpha$
for some \(c'>0\) and some \(\alpha\in\mathcal A_{n-1}\) satisfying the
required bounds.
Finally, suppose that \(c<0\) and \(n\) is even. In this case, \(q\) is negative
both on \((-\infty,r_1)\) and on \((r_n,+\infty)\). Since neither interval can
contain a point of \(K_N\), we must have $r_1=-N$ and $r_n=N.$
Therefore,
\[
\begin{aligned}
q(x)
&=
c(x+N)(x-N)\prod_{i=2}^{n-1}(x-r_i) = 
(-c)(N+x)(N-x)\prod_{i=2}^{n-1}(x-r_i).
\end{aligned}
\]
For $n=2$, this gives $q(x)=(-c)(N^2-x^2)$.
For $n\ge4$, the remaining negative intervals are bounded by consecutive
pairs of roots, so $r_{2k+1}=r_{2k}+1$ for $k=1,\ldots,(n-2)/2$.
Thus $(r_2,\ldots,r_{n-1})\in\mathcal A_{n-2}$, with
$r_2\ge -N+1$ and $r_{n-1}\le N-1$, giving the required factorization.
The four alternatives are disjoint, since they correspond to the four possible
combinations of the sign of the leading coefficient and the parity of \(n\).
Within each alternative, the ordered root tuple uniquely determines
the normalized polynomial.
\end{proof}
The four families introduced in \thref{def:four-families} give exactly
the defining inequalities of the facets of $\cS_n(K_N)$.
\begin{theorem}[Solvability criterion for the truncated $K_N$-moment problem]
\thlabel{thm:Z-finite-support}
Let $N,n\in\mathbb{N}$ with \(2N\ge n\), and \(m=(m_1,\dots,m_n)\in\mathbb{R}^n\).
Then \(m\) admits a \(K_N\)-representing probability measure if and only if
\[
L_m(q)\ge 0, \quad \forall q\in \cE_{n,N}\cup\cQ_{n,N}^+\cup\cQ_{n,N}^-\cup\cR_{n,N}.
\]
Moreover, none of these inequalities is redundant.
\end{theorem}

\begin{proof}
Necessity is immediate, since every polynomial in \(\cE_{n,N}\cup\cQ_{n,N}^+\cup\cQ_{n,N}^-\cup\cR_{n,N}\) is non-negative on
\(K_N\), and so \thref{prop:positivity-support}-\textup{(i)} applies.
Conversely, \thref{prop:finite-support-polytope} gives
\[
\cS_n(K_N)=\conv (V_N),
\mbox{where }V_N:=\{v_n(k):k\in\{-N,\dots,N\}\}
\]
and $v_n$ denotes the projected moment curve of order $n$
(see \thref{def:momentcurve}). Since $2N\ge n$, the set $V_N$
contains $2N+1\ge n+1$ points. Any $n+1$ distinct points on the
projected moment curve are affinely independent, so $\conv(V_N)$
is full-dimensional in $\R^n$. Thus, we can apply
\thref{lem:convex-hull-halfspaces} to $V_N$.
 Therefore, $\cS_n(K_N)$ is the intersection of the closed half-spaces
containing $V_N$ whose boundary hyperplanes contain $n$ points of $V_N$.
For each such half-space, we write its defining inequality as $L_m(q)\ge0$
for some $q\in\R[x]_{\le n}$, by using the identification of affine
functionals with polynomials. Then $q$ is non-negative on $K_N$ and has at least
$n$ distinct zeros there. Hence, its degree is exactly $n$ and it has
exactly $n$ distinct such zeros. Conversely, each polynomial with these properties
defines one of the stated half-spaces.
After normalization by the absolute value of the leading coefficient,
\thref{lem:KN-zero-classification} gives exactly the following four types:
\begin{itemize}
    \item $q=p$ with $p\in\cE_{n,N}$;
    \item $q=(N-x)p$ with $(N-x)p\in\cQ_{n,N}^{+}$;
    \item $q=(N+x)p$ with $(N+x)p\in\cQ_{n,N}^{-}$;
    \item $q=(N+x)(N-x)p$ with $(N+x)(N-x)p\in\cR_{n,N}$.
\end{itemize}
Thus, we obtain
\[
\cS_n(K_N)
=
\bigcap_{q\in\cE_{n,N}\cup\cQ_{n,N}^+\cup\cQ_{n,N}^-\cup\cR_{n,N}}
\{m\in\mathbb R^n:L_m(q)\ge 0\}
\]
as needed. Distinct normalized polynomials define distinct half-spaces: indeed,
equal half-spaces would give positively proportional affine functionals and, hence, positively proportional polynomials, which normalization excludes.
The non-redundancy now follows from \thref{lem:convex-hull-halfspaces}.
\end{proof}

\begin{example}
Let $K_4=\{-4,\dots,4\}$. The polytope $\cS_2(K_4)=\conv\{(k,k^2):k\in K_4\}$
has eight lower facets induced by
$(x-k)(x-k-1)\in\cE_{2,4}$ with $k\in\{-4,\dots,3\}$, and one upper
facet induced by $16-x^2\in\cR_{2,4}$. Hence,
\[
\cS_2(K_4)=
\left\{(m_1,m_2)\in\R^2:
\begin{aligned}
m_2-(2k+1)m_1+k(k+1)&\ge0,
&&k=-4,\dots,3,\\
m_2&\le16
\end{aligned}
\right\},
\]
as shown in Figure~\ref{fig:moment-polytope-K4}.

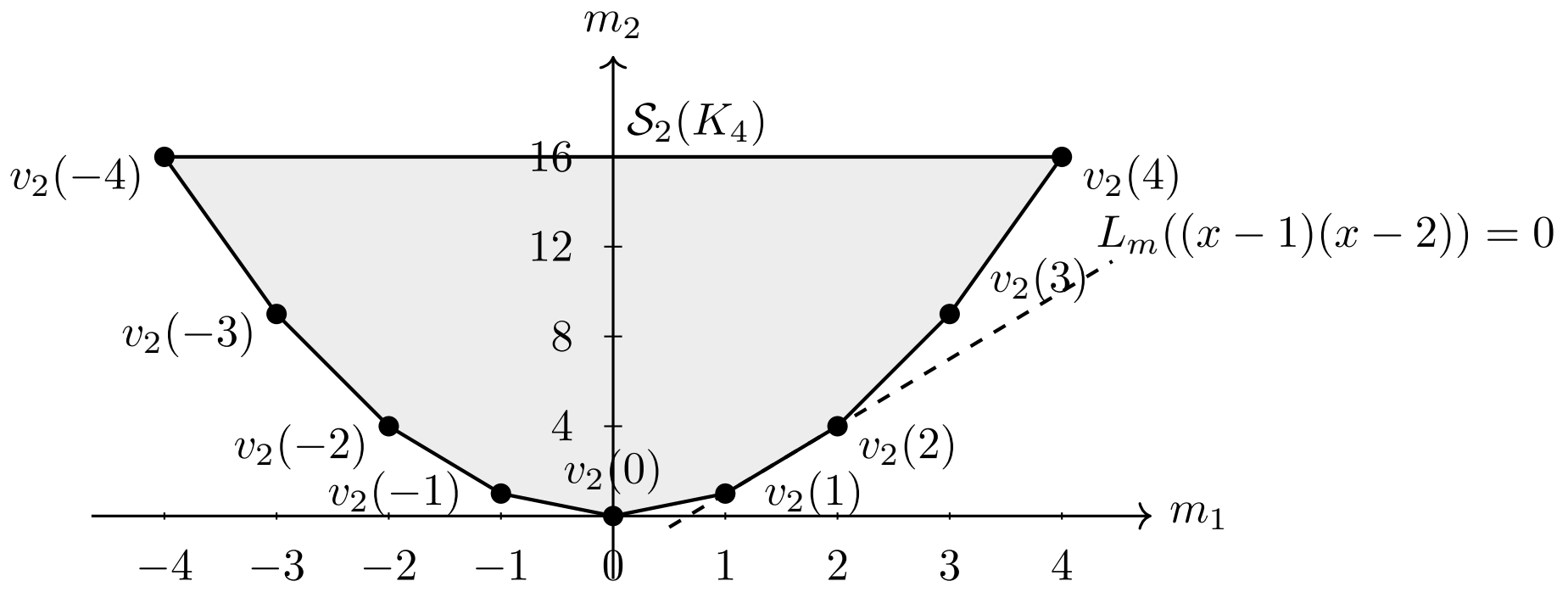
\begin{figure}[H]
\centering
\begin{tikzpicture}[
    x=1.25cm,
    y=0.25cm,
    scale = 0.7,
    every node/.style={font=\small},
    axis/.style={->, semithick},
    tick/.style={thin},
    polytope/.style={fill=black!7, draw=black, line width=0.8pt},
    vertex/.style={circle, fill=black, inner sep=1.6pt},
    intrinsic/.style={black, dashed, line width=0.8pt},
    boundary/.style={black, densely dotted, line width=0.9pt}
]

\filldraw[polytope]
    (-4,16) --
    (-3,9) --
    (-2,4) --
    (-1,1) --
    (0,0) --
    (1,1) --
    (2,4) --
    (3,9) --
    (4,16) --
    cycle;

\foreach \x/\y in {
    -4/16,
    -3/9,
    -2/4,
    -1/1,
     0/0,
     1/1,
     2/4,
     3/9,
     4/16
} {
    \node[vertex] at (\x,\y) {};
}

\node[below left=1pt] at (-4,17) {$v_2(-4)$};
\node[below left=1pt] at (-3,10)  {$v_2(-3)$};
\node[below left=1pt] at (-2,5)  {$v_2(-2)$};
\node[left =1pt] at (-1.15,1) {$v_2(-1)$};
\node[above = 1pt] at (0,0) {$v_2(0)$};
\node[right =1pt] at (1.15,1) {$v_2(1)$};
\node[below right=1pt]  at (2,5)   {$v_2(2)$};
\node[anchor=west] at (3.2,10.5) {$v_2(3)$};
\node[below right=1pt]  at (4,17)  {$v_2(4)$};

\draw[intrinsic] (0.5,-0.5) -- (4.45,11.35);
\node[anchor=west] at (4.15,12.5)
    {$L_m\!\left((x-1)(x-2)\right)=0$};

\node[align=center] at (0.75,17.5)
    {$\cS_2(K_4)$};

\draw[axis] (-4.65,0) -- (4.8,0) node[right] {$m_1$};
\draw[axis] (0,-2.8) -- (0,20.5) node[above] {$m_2$};

\foreach \x in {-4,...,4} {
    \draw[tick] (\x,0.16) -- (\x,-0.16);
    \node[below=4pt] at (\x,0) {$\x$};
}

\foreach \y in {4,8,12,16} {
    \draw[tick] (0.08,\y) -- (-0.08,\y);
    \node[left=5pt] at (0,\y) {$\y$};
}

\end{tikzpicture}
\caption{In light grey, the moment polytope $\cS_2(K_4)$ associated with
$K_4=\{-4,\dots,4\}$. The lower facets are induced by polynomials in
$\cE_{2,4}$, whereas the polynomial $(4+x)(4-x)\in\cR_{2,4}$ gives rise
to the upper facet.}
\label{fig:moment-polytope-K4}
\end{figure}
\end{example}
\section{\texorpdfstring{The truncated $\Z-$moment problem}{The truncated Z-moment problem}}
\label{sec:Z-moment-problem}
We now pass from the finite sets $K_N$ to the full integer lattice.
By \thref{prop:atomic-reduction}, a vector $m=(m_1,\ldots,m_n)$ admits
a $\Z-$representing probability measure if and only if it admits one
supported on $K_N$ for some $N\in\N$ with $2N\ge n$. Among the four families of
facet inequalities in \thref{thm:Z-finite-support}, the family $\mathcal{E}_{n,N}$ is contained in
$\mathcal{E}_n$ for every $N$, while the boundary families
$\mathcal{Q}^{\pm}_{n,N}$ and $\mathcal{R}_{n,N}$ depend on $N$. The criterion below
therefore consists of non-negativity of $L_m$ on $\mathcal{E}_n$ or
$\mathcal{E}_{n-1}$, together with conditions ensuring that the boundary
inequalities hold for some sufficiently large $N$, as obtained in Appendix~\ref{app:core-technical-results}. 

We first observe that non-negativity of the Riesz functional on the polynomials of even degree propagates to all lower even degrees.

\begin{lemma}
\thlabel{lem:positivityZ}
Let $r\in\mathbb N$ and $m=(m_1,\dots,m_{2r})\in\mathbb R^{2r}$. If $L_m(p)\geq0$ for every $p\in\cE_{2r}$, then $L_m(q)\geq0$ for every $q\in\cE_{2r-2}$.
Moreover, if for any such $q\in\cE_{2r-2}$ we have $L_m(q)=0$, then
$L_m(xq)=0$ and $L_m(x^2q)\ge0$.
\end{lemma}

\begin{proof}
Fix $q\in\cE_{2r-2}$. For every $a\in\Z$ such that $\{a,a+1\}\cap Z(q)=\varnothing$, we have $B_aq\in\cE_{2r}$ and, therefore,
\[
0
\leq
L_m(B_aq)
=
L_m(x^2q)-(2a+1)L_m(xq)+a(a+1)L_m(q).
\]
If $L_m(q)<0$, then the right-hand side tends to $-\infty$ as $a\to\infty$ along the integers, which is impossible. Hence, $L_m(q)\geq0$.
If $L_m(q)=0$, the right-hand side is affine in $a$ and has to be
non-negative for arbitrarily large positive and negative integers $a$.
This implies $L_m(xq)=0$, and the remaining constant is $L_m(x^2q)\ge0$.
\end{proof}

For any $p\in\bigcup_{r\ge0}\cE_{2r}$, we define
\begin{equation}
\label{eq:rho(P)-def}
\rho(1):=0,\qquad
\rho(p):=\max\{|t|:p(t)=0\}\quad\text{if }p\ne1.
\end{equation} 
\begin{theorem}[Solvability criterion for the truncated \(\Z-\)moment problem]
\thlabel{thm:Z-main}
\leavevmode\par\noindent
Let \(m=(m_1,\dots,m_n)\in\mathbb{R}^n\).
\begin{enumerate}[label=(\roman*),leftmargin=1.2cm]
\item If \(n\) is odd, then \(m\) admits a \(\Z-\)representing probability measure if and only if:
\begin{enumerate}[label=(\arabic*),leftmargin=1.2cm]
\item $L_m(p)\geq0$ for every $p\in\cE_{n-1};$
\item \(\forall p\in\cE_{n-1}\) such that \(L_m(p)=0\), $L_m(xp)=0$;
\item \(\forall p\in\cE_{n-3}\) such that \(L_m(p)=0\), $L_m(x^2p)=0$.
\end{enumerate}
\item If \(n\) is even, then \(m\) admits a \(\Z-\)representing probability measure if and only if:
\begin{enumerate}[label=(\arabic*),leftmargin=1.2cm]
\item $L_m(p)\geq0$ for every $p\in\cE_n;$
\item \(\forall p\in\cE_{n-2}\) such that \(L_m(p)=0\), $L_m(x^2p)=0$.
\end{enumerate}
\end{enumerate}
Conditions involving \(\cE_j\) with \(j<0\) are understood to be void.
\end{theorem}

\begin{proof}

\emph{Necessity.}
Let $\mu$ be a $\Z-$representing probability measure for $m$.
Since the polynomials in the families $\cE_j$ are non-negative
on $\Z$, conditions \textup{(1)} in \textup{(i)} and \textup{(ii)} follow from
\thref{prop:positivity-support}. If such a polynomial $p$ satisfies
$L_m(p)=0$, the same proposition gives $\supp\mu\subseteq Z_{\Z}(p)$.
Consequently $L_m(x^jp)=0$ for every integer $j\ge0$ with
$\deg p+j\le n$, which proves \textup{(2)}--\textup{(3)} in \textup{(i)} and \textup{(2)} in \textup{(ii)}.

\emph{Sufficiency for odd $n$.}
Let us assume conditions \textup{(1)}--\textup{(3)} in \textup{(i)}.
For $n=1$, we choose $a\in\Z$ with $a\le m_1\le a+1$ and set
$\lambda:=m_1-a$. Then
$\mu:=(1-\lambda)\delta_a+\lambda\delta_{a+1}$ represents $m_1$ on $\Z$.
Suppose henceforth that $n\ge3$.
Since $\cE_{n,N}=\cR_{n,N}=\varnothing$, it is enough, by
\thref{thm:Z-finite-support}, to find an integer $N$ with $2N\ge n$
for which $L_m$ is non-negative on
$\cQ_{n,N}^{+}\cup\cQ_{n,N}^{-}$.

Let us write $n=2r+1$ for some $r\geq 1$. Condition \textup{(1)} gives
$L_m(p)\ge0$ for every $p\in\cE_{2r}$, while condition
\textup{(3)} states that
\[
L_m(H)=0\quad\Longrightarrow\quad L_m(x^2H)=0,
\qquad H\in\cE_{2r-2}.
\]
These are exactly the hypotheses of
\thref{thm:thm3.2-gen}\textup{(b)} with index $r$.
We therefore obtain a finite set
$\mathcal F:=\widetilde{\cE}_{n-1}^{m}$, independent of $N$
and of the sign of the boundary factor. More precisely, for
every $p\in\cE_{n-1}$ there is a $\widetilde p\in\mathcal F$
with $\rho(\widetilde p)\le\rho(p)$ such that
\[
L_m((N+\varepsilon x)\widetilde p)
\le L_m((N+\varepsilon x)p)
\]
for every integer $N\ge\rho(p)$ and every
$\varepsilon\in\{-1,1\}$. Thus, the representative can be
chosen before either $N$ or the sign is fixed.
For each $\widetilde p\in\mathcal F$, condition \textup{(1)} gives
$L_m(\widetilde p)\ge0$. If this value is zero, condition
\textup{(2)} gives $L_m(x\widetilde p)=0$ as well.
We define the finite constant
\[
C_1:=\max\left(\{0\}\cup
\left\{\frac{|L_m(xp)|}{L_m(p)}:
p\in\mathcal F,\ L_m(p)>0\right\}\right)
\]
and we choose
\[
N:=\left\lceil\max\left\{\frac n2,
\max_{p\in\mathcal F}\rho(p),\ C_1\right\}\right\rceil.
\]
Then $2N\ge n$ and, for every $\widetilde p\in\mathcal F$ we have
\[
N\ge\rho(\widetilde p),\qquad
N L_m(\widetilde p)\ge |L_m(x\widetilde p)|.
\]
Every $q\in\cQ_{n,N}^{+}\cup\cQ_{n,N}^{-}$ has the form
$q=(N+\varepsilon x)p$ for some $p\in\cE_{n-1}$ with $\rho(p)\le N$.
Using a representative $\widetilde{p}\in\mathcal{F}$ we obtain
\[
0\le N L_m(\widetilde p)+\varepsilon L_m(x\widetilde p)
=L_m((N+\varepsilon x)\widetilde p)
\le L_m(q).
\]
Hence, by \thref{thm:Z-finite-support} there is a $K_N$-representing
probability measure for $m$. In particular this is a $\Z-$representing probability measure for $m$.

\emph{Sufficiency for even $n$.}
Let us assume conditions \textup{(1)}--\textup{(2)} in \textup{(ii)}.
Condition \textup{(1)} already gives $L_m\ge0$ on
$\cE_{n,N}\subseteq\cE_n$ for every $N$ with $2N\ge n$.
Since $\cQ_{n,N}^{+}=\cQ_{n,N}^{-}=\varnothing$, only the
inequalities associated with $\cR_{n,N}$ remain to be verified.
Moreover, by \thref{lem:positivityZ} we have $L_m(p)\ge0$ for every
$p\in\cE_{n-2}$.
For $n=2$, no finite-reduction theorem is needed: $\cE_0=\{1\}$,
and the only remaining inequality is $N^2\ge m_2$.
We choose an integer $N\ge1$ satisfying this inequality and we apply
\thref{thm:Z-finite-support}. We may therefore assume $n\ge4$.

Condition \textup{(1)} allows us to apply
\thref{thm:thm3.2-gen}\textup{(c)}. Thus, there is a finite set
$\widetilde{\cE}_{n-2}^{m}\subseteq\cE_{n-2}$, independent of $N$,
such that every $p\in\cE_{n-2}$ has a representative
$\widetilde p\in\widetilde{\cE}_{n-2}^{m}$ satisfying
$\rho(\widetilde p)\le\rho(p)$ and
\[
L_m((N^2-x^2)\widetilde p)\le L_m((N^2-x^2)p)
\]
for every integer $N\ge\rho(p)$.
Consider $\mathcal F:=\widetilde{\cE}_{n-2}^{m}$.
For every $\widetilde p\in\mathcal F$ we have $L_m(\widetilde p)\ge0$ and,
if this quantity vanishes, condition \textup{(2)} forces
$L_m(x^2\widetilde p)=0$. Thus, by setting
\[
C_2:=\max\left(\{0\}\cup
\left\{\frac{L_m(x^2p)}{L_m(p)}:
p\in\mathcal F,\ L_m(p)>0\right\}\right),
\]
we may choose
\[
N:=\left\lceil\max\left\{\frac n2,
\max_{p\in\mathcal F}\rho(p),\ \sqrt{C_2}\right\}\right\rceil.
\]
This ensures $2N\ge n$ and, for every $\widetilde p\in\mathcal F$ we have
\[
N\ge\rho(\widetilde p),\qquad
N^2 L_m(\widetilde p)\ge L_m(x^2\widetilde p).
\]
Now let $R\in\cR_{n,N}$. By \thref{def:four-families} we have
$R=(N^2-x^2)p$ for some $p\in\cE_{n-2}$ whose roots lie strictly
between $-N$ and $N$, so that $\rho(p)<N$.
There exists a representative $\widetilde{p}\in\mathcal{F}$ such that
\[
0\le N^2 L_m(\widetilde p)-L_m(x^2\widetilde p)
 =L_m((N^2-x^2)\widetilde p)
 \le L_m(R).
\]
Thus, all the inequalities in \thref{thm:Z-finite-support} hold,
and $m$ admits a $K_N$-representing probability measure.
In particular, in either parity, this is also a $\Z-$representing probability measure.
\end{proof}

\begin{example}
We now illustrate why each compatibility condition is needed.
\thlabel{ex:Z-compatibility}

\begin{itemize}

\item[$n=3$.]
Let $m=(0,0,1)$. Since $L_m(B_a)=a(a+1)\ge0$ for every $a\in\Z,$
condition \textup{(i)(1)} holds. However,
\[
L_m(B_0)=0,
\qquad
L_m(xB_0)=1,
\]
so condition \textup{(i)(2)} fails. Indeed, $m_2=0$ would force
$\mu=\delta_0$, which contradicts $m_3=1$.

\item[$n=4$.]
Let $m=(0,0,0,1)$. For every $p\in\cE_4$ $L_m(p)=1+p(0)\ge1,$
so condition \textup{(ii)(1)} even holds strictly. Nevertheless,
\[
L_m(B_0)=0,
\qquad
L_m(x^2B_0)=1
\]
and, hence, condition \textup{(ii)(2)} fails. Again, $m_2=0$ would force
$\mu=\delta_0$, which contradicts $m_4=1$.

\item[$n=5$.]
Let $m=(0,0,0,1,0)$. Conditions \textup{(i)(1)} and
\textup{(i)(2)} hold, the latter vacuously, whereas
\[
L_m(B_0)=0,
\qquad
L_m(x^2B_0)=1.
\]
Thus, condition \textup{(i)(3)} fails.

\end{itemize}

Therefore, conditions \textup{(i)(2)}, \textup{(i)(3)}, and
\textup{(ii)(2)} in \thref{thm:Z-main} cannot be omitted.
\end{example}

\section{\texorpdfstring{From extremal completions to a finite reduction}{From extremal completions to a finite reduction}}

\label{sec:extremal-completions-Z}

The criterion in \thref{thm:Z-main} involves infinitely many
polynomials. We now turn it into a finite recursive procedure, by adapting
the strategy of \cite[Section~6.1]{InfusinoKunaLebowitzSpeer2017}
to $\Z$. To this aim we will first introduce a classification of the moment vectors in 
I-representable and B-representable (see \thref{def:Z-interior-boundary}) as well as the concept of extremal $\Z-$ and Hamburger completions (see \thref{def:extremal-Z-completion} and \thref{def:Hamburger-completion}).

\subsection{Extremal $\Z-$completion and induction}

Let us denote $m^{(k)}:=(m_1,\ldots,m_k)$ and $d_k:=2\lfloor k/2\rfloor$.
We use the following terminology, analogous to the distinction for
$K = \mathbb N_0$ in \cite{InfusinoKunaLebowitzSpeer2017}.

\begin{definition}
\thlabel{def:Z-interior-boundary}
A vector $m^{(k)}\in\cS_k(\Z)$ is \emph{I-representable} on $\Z$ if
$L_{m^{(k)}}(p)>0$ for every $p\in\cE_{d_k}$ and
\emph{B-representable} on $\Z$ otherwise. Thus, a B-representable vector
is associated with at least one polynomial $p\in\cE_{d_k}$ such that $L_{m^{(k)}}(p)=0$.
\end{definition}

Strict positivity also descends to all lower even degrees:
if $m^{(k)}$ is I-representable on $\Z$, then
\begin{equation}
\label{eq:Z-lower-strict-positivity}
L_{m^{(k)}}(q)>0,
\qquad q\in\bigcup_{0\le 2j\le k}\cE_{2j}.
\end{equation}
Indeed, a polynomial $q$ with zero Riesz value would force every $\Z-$representing probability measure to be
supported on $Z(q)$. After multiplying $q$ by sufficiently distant, pairwise
disjoint lattice blocks we would obtain a polynomial $p\in\cE_{d_k}$ with zero Riesz value, contradicting the I-representability.

\begin{definition}
\thlabel{def:Z-minimizing-polynomial}
    For an odd-order vector $m:=m^{(2r-1)}$ and a polynomial $p\in\cE_{2r}$, we set
$$
\Phi_m(p):=L_{(m_1,\ldots,m_{2r-1},0)}(p).
$$
We call a minimizer of $\Phi_m$ among the polynomials in $\cE_{2r}$ a
\emph{minimizing polynomial of order $2r$}.
\end{definition}
Note that, since every polynomial in $\cE_{2r}$ is monic, $L_{(m_1,\dots,m_{2r-1},t)}(p)=t+\Phi_m(p)$ for every $t\in\R$. This implies that for every $t\in\mathbb{R}$
\[
\argmin_{p\in\cE_{2r}} \Phi_m(p) = \argmin_{p\in\cE_{2r}} L_{(m_1,\dots,m_{2r-1},t)}(p).
\]

\begin{proposition}[Inductive solvability criterion]
\thlabel{prop:Z-inductive}
Every $m^{(1)}\in\R$ is I-representable on $\Z$. For $k\ge2$, a non-$\Z-$representable
$m^{(k-1)}$ has no $\Z-$representable extension. If $m^{(k-1)}$ is $\Z-$representable,
the following alternatives hold.
\begin{enumerate}[label=\textup{(\roman*)},leftmargin=1cm]
\item If $m^{(k-1)}$ is I-representable on $\Z$ and $k$ is odd, every
$m_k\in\R$ gives an I-representable extension.
\item If $m:=m^{(2r-1)}$ is I-representable on $\Z$, there exists at least one minimizing polynomial $p_{2r}^{(m)}\in\cE_{2r}$ of order $2r$. We set
\[
\tau_{\Z}(m):=-\Phi_m\bigl(p_{2r}^{(m)}\bigr),
\]
and the extension $(m_1,\dots,m_{2r-1},t)$ is $\Z-$representable if and only if $t\ge\tau_{\Z}(m)$.
It is I-representable if $t>\tau_{\Z}(m)$ and B-representable if
$t=\tau_{\Z}(m)$. The extremal $\Z-$extension has a unique $\Z-$representing probability
measure, denoted by $\mu_{\Z}(m)$, with
$\supp\mu_{\Z}(m)\subseteq Z(p_{2r}^{(m)})$ for every minimizing
polynomial $p_{2r}^{(m)}\in\mathcal{E}_{2r}$ of order $2r$.
\item If $m^{(k-1)}$ is B-representable on $\Z$, its $\Z-$representing probability
measure $\mu$ is unique, and the extension is $\Z-$representable if and only if
\begin{equation}
\label{eq:Z-boundary-extension}
m_k=\int_{\Z}x^k\,d\mu(x).
\end{equation}
When this holds, the extension is again B-representable. Equivalently,
if $p\in\cE_{d_{k-1}}$ satisfies $L_{m^{(k-1)}}(p)=0$, the condition to be $\Z-$representable is given by
$L_{m^{(k)}}(x^{k-d_{k-1}}p)=0$.
\end{enumerate}
\end{proposition}

\begin{proof}
Let $a:=\lfloor m_1\rfloor$ and $\theta:=m_1-a\in[0,1)$. The probability
measure $\mu=(1-\theta)\delta_a+\theta\delta_{a+1}$ is supported on $\Z$
and has first moment $m_1$. Since $\cE_0=\{1\}$, every $m^{(1)}$ is
I-representable. Moreover, a representing measure for a moment vector
also represents all its prefixes. Thus, a non-representable prefix cannot
have a $\Z-$representable extension.

\emph{(i).} Since $k$ is odd, we have
$d_k=d_{k-1}=k-1$. Hence, for every
$p\in\cE_{d_k}$
\[
L_{m^{(k)}}(p)
=L_{m^{(k-1)}}(p)>0.
\]
By \eqref{eq:Z-lower-strict-positivity}, the zero-propagation
conditions in \thref{thm:Z-main} are void. Hence, the extension is
$\Z-$representable and therefore I-representable, independently of $m_k$.

\emph{(ii).} We first show that $\Phi_m$ attains its minimum on
$\cE_{2r}$. For $r\ge2$, this follows from
\eqref{eq:Z-lower-strict-positivity} and
\thref{thm:thm3.2-gen}\textup{(a)}. For $r=1$, let
$a:=\lfloor m_1\rfloor$ and $\theta:=m_1-a$. Then
\[
\Phi_m(B_{a+h})-\Phi_m(B_a)=h(h+1-2\theta)\ge0,
\qquad h\in\Z,
\]
so $B_a$ is a minimizer.

For every $p\in\cE_{2r}$ and every $t\in\R$, monicity gives
\[
L_{(m,t)}(p)=t+\Phi_m(p),
\qquad
\min_{p\in\cE_{2r}}L_{(m,t)}(p)
=t+\Phi_m\bigl(p_{2r}^{(m)}\bigr)
=t-\tau_{\Z}(m).
\]
Since $m$ is I-representable, the zero-propagation condition in the
even-order part of \thref{thm:Z-main} is void. Thus, the extension is
$\Z-$representable exactly when $t\ge\tau_{\Z}(m)$, and strict
positivity on $\cE_{2r}$ holds exactly when $t>\tau_{\Z}(m)$.

At equality, $L_{(m,t)}(p_{2r}^{(m)})=0$ for every minimizing polynomial $p_{2r}^{(m)}$. By
\thref{prop:positivity-support}, every representing measure is supported
on the $2r$ distinct zeros of $p_{2r}^{(m)}$. The
weights at these points are uniquely determined by the moments of orders $0,\ldots,2r-1$, through an
invertible Vandermonde system, so the representing
measure is unique.

\emph{(iii).} Set $d:=d_{k-1}$ and choose $p\in\cE_d$ such that
$L_{m^{(k-1)}}(p)=0$. Since $p$ is non-negative on $\Z$,
\thref{prop:positivity-support} forces every representing measure of
$m^{(k-1)}$ to be supported on $Z(p)$. This set consists of $d\le k-1$
distinct points, so the moments of orders $0,\ldots,d-1$ uniquely
determine the weights. In particular, the representing measure $\mu$
is unique. Consequently, the extension is $\Z-$representable if and
only if $m_k=\int_{\Z}x^k\,d\mu(x)$.

To obtain the equivalent condition, write
\[
x^{k-d}p=x^k+R,\qquad \deg R\le k-1.
\]
Since $\mu$ is supported on $Z(p)$,
\[
0=\int_{\Z}x^{k-d}p(x)\,d\mu(x)
=\int_{\Z}x^k\,d\mu(x)+L_{m^{(k-1)}}(R).
\]
Therefore, $L_{m^{(k)}}(x^{k-d}p)=0$ is equivalent to
\eqref{eq:Z-boundary-extension}. Finally, the extension remains
B-representable: if $k$ is odd, $p\in\cE_{d_k}$ still has zero Riesz
value; if $k$ is even, the same holds for $pB_a\in\cE_{d_k}$, where
$B_a$ is any lattice block with zeros disjoint from $Z(p)$.
\end{proof}
\begin{definition}[Extremal $\Z-$completion]
\thlabel{def:extremal-Z-completion}
Let $m=m^{(2r-1)}$ be I-representable on $\Z$.
The vector
$
(m,\tau_{\Z}(m))
$
is called the \emph{extremal $\Z-$completion} of $m$,
and its unique $\Z-$representing probability measure
$\mu_{\Z}(m)$ is called the
\emph{extremal $\Z-$representing measure}.
\end{definition}
\begin{example}
    In particular,
for $r=1$ the above proof gives
\begin{equation}
\label{eq:Z-order-two}
\begin{aligned}
p_2^{(m)}&=B_a,&
\tau_{\Z}(m)&=(2a+1)m_1-a(a+1)
              =m_1^2+\theta(1-\theta),\\
\mu_{\Z}(m)&=(1-\theta)\delta_a+\theta\delta_{a+1},&
a&=\lfloor m_1\rfloor,\quad \theta=m_1-a.
\end{aligned}
\end{equation}
When $m_1$ is
an integer, then $a = m_1$ and the extremal $\Z-$representing measure is $\delta_a$. In this case both
$B_{a-1}$ and $B_a$ are minimizing polynomials.
\end{example}

\subsection{Hamburger completion and interlacing}\label{subsec:hamburger-completion}
The following elementary completion rule for lattice root sets will be
used both for the Hamburger construction and for the recursive step.

\begin{lemma}[Compatible completion of lattice root sets]
\thlabel{lem:Z-compatible-completion}
Let $s\in\N$. A union of at most $s$ sets of the form $\{a,a+1\}$, with $a\in\Z$, is contained in the zero set of some polynomial in $\cE_{2s}$. Consequently:
\begin{enumerate}[label=\textup{(\roman*)},leftmargin=1cm]
\item every $S\subseteq\Z$ with $|S|\le s$ is contained in such a zero set;
\item for every $r\ge1$, if $S\subseteq Z(R)$ for some
$R\in\cE_{2r-2}$ and $a\in\Z$, then
$S\cup\{a,a+1\}\subseteq Z(q)$ for some
$q\in\cE_{2r}$.
\end{enumerate}
The completion can be constructed by a finite procedure.
\end{lemma}

\begin{proof}
Let $U=\bigcup_{j=1}^{k}\{a_j,a_j+1\}$ with $k\le s$.
We construct a collection of disjoint adjacent pairs
covering $U$ by the following simple procedure.

Take the smallest element $b$ of $U$ not yet covered
and add the pair $\{b,b+1\}$. Repeat until every
element of $U$ is covered. By construction, the
selected pairs are disjoint.

We claim that this procedure uses at most $k$ pairs.
Indeed, the starting points of the selected pairs
differ by at least two. Each of them belongs to
one of the original $k$ pairs, and no original pair
can contain two such starting points. Hence, the
number of selected pairs cannot exceed $k\le s$.

If fewer than $s$ pairs have been selected, we add
further disjoint adjacent pairs sufficiently far
to the right until we have exactly $s$ pairs,
say $\{c_j,c_j+1\}$ with $j\in\{1,\ldots,s\}$. Then
\[
q(x)=\prod_{j=1}^{s}(x-c_j)(x-c_j-1)
\in\cE_{2s},
\]
and $U\subseteq Z(q)$. This also gives the
asserted finite construction.

For \textup{(i)}, each $z\in S$ belongs to the
pair $\{z,z+1\}$. Since $|S|\le s$, the first
part applies.

For \textup{(ii)}, the zero set of
$R\in\cE_{2r-2}$ is the union of $r-1$
disjoint adjacent pairs. Adding $\{a,a+1\}$
gives a union of at most $r$ adjacent pairs.
Applying the first part with $s=r$ yields
$q\in\cE_{2r}$ such that $S\cup\{a,a+1\}\subseteq Z(q).$
\end{proof}

\begin{example}
Let $s=4$ and consider
\[
U=\{0,1\}\cup\{1,2\}\cup\{5,6\}.
\]
The procedure described in \thref{lem:Z-compatible-completion} first selects $\{0,1\}$,
then $\{2,3\}$ to cover the remaining point
$2$, and finally $\{5,6\}$. Since we have
obtained only three disjoint pairs, we add
$\{7,8\}$ to reach $s=4$. Thus, $U\subseteq
\{0,1\}\cup\{2,3\}\cup
\{5,6\}\cup\{7,8\},$
and the corresponding polynomial is
\[
q=B_0B_2B_5B_7\in\cE_8,
\qquad U\subseteq Z(q).
\]
\end{example}

Let us fix an I-representable vector $m:=m^{(2r-1)}$, and allow the
representing measure to be supported on $\R$. This gives rise to the extremal Hamburger
completion, which can be computed explicitly (see \cite{CurtoFialkow1991}). To characterize the latter completion, we introduce
\[
A_{r-1}:=(m_{i+j})_{i,j=0}^{r-1},
\qquad h_r:=(m_r,\ldots,m_{2r-1})^T.
\]
\begin{definition}[Extremal Hamburger completion]
\thlabel{def:Hamburger-completion}
Let $m=m^{(2r-1)}$ be I-representable on $\Z$, and set
\[
\tau_{\R}(m):=\inf\{t\in\R:(m,t)\in\cS_{2r}(\R)\}.
\]
The vector $(m,\tau_{\R}(m))$
is called the \emph{extremal Hamburger completion} of $m$.
\end{definition}
The following proposition characterizes the extremal
Hamburger completion and relates it to the extremal
$\Z-$completion.
\begin{proposition}[Characterization of the extremal Hamburger completion]
\thlabel{prop:Z-hamburger-completion}
If $m = (m_1,\dots,m_{2r-1})$ is I-representable on $\Z$, then the matrix $A_{r-1}$ is positive definite. Moreover, set
\begin{equation}
\label{eq:Z-hamburger-system}
q_r(x):=x^r+\sum_{i=0}^{r-1}c_ix^i,
\qquad
(c_0,\ldots,c_{r-1})^T:=-A_{r-1}^{-1}h_r.
\end{equation}
Then 
\begin{equation}
\label{eq:tau_R}
    \tau_{\R}(m)=h_r^TA_{r-1}^{-1}h_r,
\end{equation}
and $(m,t)\in\cS_{2r}(\R)$ if and only if $t\ge\tau_{\R}(m)$.
At equality, the unique $\R-$representing probability measure $\mu_\mathbb{R}(m)$ has
exactly $r$ atoms, corresponding to the distinct real zeros $y_1<\cdots<y_r$ of $q_r$, and its positive weights are uniquely determined by the moments
of orders $0,\ldots,r-1$.
Furthermore, $\tau_{\R}(m)\le\tau_{\Z}(m)$ and, if every $y_j$ is an integer, then
$\tau_{\R}(m)=\tau_{\Z}(m)$ and $\mu_\mathbb{R}(m)=\mu_{\Z}(m)$.
\end{proposition}

\begin{proof}
Since $m$ is $\Z-$representable, $A_{r-1}$ is positive
semidefinite. If it were singular, there would exist
a nonzero polynomial $R$ of degree at most $r-1$
such that
\[
L_m(R^2)=0.
\]
Thus, every $\Z-$representing measure of $m$ would
be supported on at most $r-1$ points. By
\thref{lem:Z-compatible-completion}\textup{(i)},
for $r\ge2$ there would exist $p\in\cE_{2r-2}$
with $L_m(p)=0$, which contradicts I-representability, while for $r=1$, $A_0=(1)$. Hence, $A_{r-1}$ is positive
definite. By \cite[Theorem~3.8(ii)]{CurtoFialkow1991}, we have $\tau_{\R}(m)=h_r^TA_{r-1}^{-1}h_r,$
and $(m,t)\in\cS_{2r}(\R)$ if and only if
$t\ge\tau_{\R}(m)$. Furthermore, the extremal
Hamburger completion admits a unique $\R-$representing
probability measure consisting of exactly $r$ atoms. By \eqref{eq:Z-hamburger-system} and \eqref{eq:tau_R}, we obtain
\begin{equation}
\label{eq:Z-hamburger-square}
L_{(m,t)}(q_r^2)
=t-h_r^TA_{r-1}^{-1}h_r
=t-\tau_{\R}(m).
\end{equation}
If $t = \tau_{\R}(m)$, then the support of
$\mu_\mathbb{R}(m)$ has to be contained in $Z(q_r)$.
Since the measure has exactly $r$ atoms,
these are precisely the $r$ distinct real
zeros of $q_r$. Their weights are uniquely
determined by the Vandermonde system.

Finally, $\tau_{\R}(m)\le\tau_{\Z}(m)$, since
every $\Z-$representing measure is also an
$\R-$representing measure. If all the zeros
of $q_r$ are integers, $\mu_\mathbb{R}(m)$ is also
$\Z-$representing and, hence,
$\tau_{\R}(m)=\tau_{\Z}(m)$. The uniqueness
of the extremal Hamburger measure then also gives
$\mu_\mathbb{R}(m)=\mu_{\Z}(m)$.
\end{proof}

\begin{remark}
\thlabel{rem:Z-completion-gap}
By applying \eqref{eq:Z-hamburger-square} to
$t=\tau_{\Z}(m)$ we get
\[
\tau_{\Z}(m)-\tau_{\R}(m) = L_{(m,\tau_{\Z}(m))}(q_r^2)
=\int_{\Z}q_r(x)^2\,d\mu_{\Z}(m)(x).
\]
Thus, the increase in the minimum moment caused by restricting
the support to $\Z$ is measured by the integral of $q_r^2$
against the extremal $\Z-$representing measure. In particular,
\[
\tau_{\Z}(m)=\tau_{\R}(m)
\quad\Longleftrightarrow\quad
\{y_1,\ldots,y_r\}\subseteq\Z.
\]
Indeed, if the two extremal values coincide, $\mu_{\Z}(m)$ also represents
the extremal Hamburger completion. Its uniqueness guarantees
$\mu_{\Z}(m)=\mu_\mathbb{R}(m)$, and so all the nodes $y_j$ are integers.
The converse is the last assertion of
\thref{prop:Z-hamburger-completion}.
\end{remark}

The interlacing argument below is the two-sided analogue of
\cite[Proposition~6.3 and Corollary~6.4]{InfusinoKunaLebowitzSpeer2017},
with the
truncated Stieltjes moment problem replaced by the truncated Hamburger moment problem. It compares the supports of the extremal
$\Z-$representing and Hamburger measures.

\begin{proposition}[Hamburger--$\Z$ interlacing]
\thlabel{prop:Z-hamburger-interlacing}
Let $m=(m_1,\ldots,m_{2r-1})$ be I-representable on $\Z$.
If $\supp\mu_\mathbb{R}(m) = \{y_1,\dots,y_r\} \nsubseteq\Z$, then there exist
$w_0,\ldots,w_r\in\supp\mu_{\Z}(m)$ such that
\begin{equation}
\label{eq:Z-interlacing}
w_0<y_1<w_1<y_2<\cdots<y_r<w_r.
\end{equation}
\end{proposition}

\begin{proof}
For brevity, let us fix $\mu:=\mu_{\Z}(m)$, $\nu:=\mu_\mathbb{R}(m)$ and
$q(x):=\prod_{i=1}^r(x-y_i)$. The two measures have the same integrals
against every polynomial of degree at most $2r-1$.
If $\supp\mu\subseteq\{y_1,\ldots,y_r\}$,
the Vandermonde system for the moments of orders $0,\ldots,r-1$ would
force their weights to coincide, which contradicts $\supp\nu\nsubseteq\Z$.
Thus, there exists a point $t\in\supp\mu\setminus\{y_1,\ldots,y_r\}$ for which, in particular, we have $\mu(\{t\})>0$. 
Now consider the polynomials
\[
F_0(x):=\frac{q(x)^2}{x-y_1},
\qquad F_r(x):=-\frac{q(x)^2}{x-y_r},
\qquad
F_j(x):=\frac{q(x)^2}{(x-y_j)(x-y_{j+1})}
\quad(0<j<r).
\]
Each denominator cancels only one copy of the corresponding factors
of $q^2$ and, hence, $F_j(y_i) = 0$ for every $j\in\{0,\dots, r\}$ and $i\in\{1,\dots, r\}$ and the polynomials $F_j$ have degree at most $2r-1$. Since the first $2r-1$ moments are the same for both measures, we have
\begin{equation}
\label{eq:mom-agreement}
    \int F_j\,d\mu=\int F_j\,d\nu=0,
\qquad \forall j\in\{0,\ldots,r\}.
\end{equation}
The polynomials $F_0$, $F_j$ ($0< j<r$), and $F_r$ are strictly
negative on
\begin{equation}
    \label{eq:intervals}
    (-\infty,y_1),\qquad (y_j,y_{j+1}),\qquad (y_r,\infty),
\end{equation}
respectively. Outside its corresponding interval, each polynomial
is non-negative and vanishes only at $y_1,\ldots,y_r$.

Suppose by contradiction that one of these intervals contains no point of $\supp\mu$, and
let $F_j$ be the corresponding polynomial. Then $F_j\ge0$ on
$\supp\mu$. The same fixed point $t$ lies outside that interval and
is not one of the nodes, which implies $F_j(t)>0$. By \eqref{eq:mom-agreement} and the latter strict inequality we get
\[
0=\int F_j\,d\mu\ge F_j(t)\mu(\{t\})>0,
\]
which is a contradiction. Thus, every interval in \eqref{eq:intervals} contains a point of $\supp{\mu}$.
We can now choose
\[
\begin{aligned}
w_0&\in\supp\mu\cap(-\infty,y_1),\\
w_j&\in\supp\mu\cap(y_j,y_{j+1}) \quad (0< j<r),\\
w_r&\in\supp\mu\cap(y_r,\infty).
\end{aligned}
\]
This choice gives \eqref{eq:Z-interlacing} and, therefore, concludes the proof.
\end{proof}

\begin{lemma}[Forced lattice block]
\thlabel{lem:Z-forced-block}
Under the hypothesis of \thref{prop:Z-hamburger-interlacing}, let us consider
$w_0,\ldots,w_r\in\supp\mu_{\Z}(m)$ satisfying \eqref{eq:Z-interlacing}.
For every minimizing polynomial
$p_{2r}^{(m)}=\prod_{i=1}^{2r}(x-\alpha_i)\in\cE_{2r}$ with
$\alpha_1<\cdots<\alpha_{2r}$, for at least one $j\in\{1,\ldots,r\}$ we have
\begin{equation}
\label{eq:Z-forced-block}
\alpha_{2j-1}=\lfloor y_j\rfloor=w_{j-1},
\qquad
\alpha_{2j}=\lfloor y_j\rfloor+1=w_j.
\end{equation}
In particular $y_j\notin\Z$, both endpoints have positive mass under
$\mu_{\Z}(m)$, and
$p_{2r}^{(m)}=B_{\lfloor y_j\rfloor}R$ for some $R\in\cE_{2r-2}$.
\end{lemma}

\begin{proof}
Let us fix any minimizing polynomial $p_{2r}^{(m)}\in\cE_{2r}$,
and write its distinct integer zeros as
$\alpha_1<\cdots<\alpha_{2r}$. By
\thref{prop:Z-inductive}\textup{(ii)},
\[
\supp\mu_{\Z}(m)\subseteq Z(p_{2r}^{(m)})
=\{\alpha_1,\ldots,\alpha_{2r}\}.
\]

Choose the largest $j\in\{1,\ldots,r\}$ such that
$\alpha_{2j-1}\le w_{j-1}$ (such an index exists because
$\alpha_1\le w_0$). If $j<r$, maximality gives
$w_j<\alpha_{2j+1}$; for $j=r$, we set
$\alpha_{2r+1}:=+\infty$. Hence, we have
\[
\alpha_{2j-1}\le w_{j-1}<w_j<\alpha_{2j+1}.
\]
The only zeros of $p_{2r}^{(m)}$ in $[\alpha_{2j-1},\alpha_{2j+1})$ are
$\alpha_{2j-1}$ and $\alpha_{2j}$. Since $w_{j-1}$ and $w_j$ are
two distinct zeros in this interval, we obtain
\[
w_{j-1}=\alpha_{2j-1},\qquad
w_j=\alpha_{2j}=\alpha_{2j-1}+1,
\]
where the last equality follows from $p_{2r}^{(m)}\in\cE_{2r}$ and the structure of the zeros of polynomials in $\cE_{2r}$.
Together with \eqref{eq:Z-interlacing}, this gives
$\alpha_{2j-1}<y_j<\alpha_{2j-1}+1$, which proves
\eqref{eq:Z-forced-block} and $y_j\notin\Z$.
Both $\alpha_{2j-1} = \lfloor y_j \rfloor$ and $\alpha_{2j} = \lfloor y_j \rfloor + 1$ have positive mass because they belong to the support
of a measure on $\Z$. After removing the factor $B_{\lfloor y_j \rfloor}$ from $p_{2r}^{(m)}$ we obtain a polynomial
$R\in\cE_{2r-2}$.
\end{proof}

Accordingly, when not all the atoms of the extremal Hamburger probability measure are integers, it is sufficient to
consider the finite set
\begin{equation}
\label{eq:Z-admissible-blocks}
\mathcal A_H(m):=\{\lfloor y_j\rfloor:1\le j\le r,\ y_j\notin\Z\}.
\end{equation}
Note that integer nodes need not be included:
the lattice block guaranteed by \eqref{eq:Z-forced-block} surrounds a non-integer
node. This implies at most $r$ branches.

\subsection{Lattice block reduction and finite minimization}

\begin{definition}[Reduction by a lattice block]
\thlabel{def:Z-block-reduction}
Let $n\ge3$ and $m=(m_1,\ldots,m_n)$ be
I-representable on $\Z$, with $m_0:=1$.
In particular, for every $a\in\Z$, by
\eqref{eq:Z-lower-strict-positivity} we have
\[
L_m(B_a)=m_2-(2a+1)m_1+a(a+1)>0.
\]
We define the reduced moment vector
$M_a(m):=(M_{a,1},\ldots,M_{a,n-2})$ by setting
$M_{a,0}:=1$ and
$$
M_{a,\ell}:=
\frac{L_m(x^\ell B_a)}{L_m(B_a)}
=\frac{m_{\ell+2}-(2a+1)m_{\ell+1}
+a(a+1)m_\ell}{L_m(B_a)},
\qquad 1\le\ell\le n-2.
$$
For any extension $(m,t)$ the same formula
with $m_{n+1}:=t$, defines the reduced vector
$M_a(m,t)$ of order $n-1$, whose first $n-2$
moments coincide with those of $M_a(m)$.
\end{definition}

\begin{lemma}[Reduction by a lattice block]
\thlabel{lem:Z-algorithm-block-reduction}
Let $n\ge3$ and let $m=m^{(n)}$ be I-representable on $\Z$.
For every $a\in\Z$, the following hold.
\begin{enumerate}[label=\textup{(\roman*)},leftmargin=1cm]
\item If $\sigma$ is a $\Z-$representing probability measure
for $m$, then
\begin{equation}
\label{eq:Z-weighted-measure}
\sigma_a(\{z\})
:=\frac{B_a(z)}{L_m(B_a)}\sigma(\{z\}),
\qquad z\in\Z,
\end{equation}
is a probability measure on $\Z$ whose moment vector is $M_a(m)$.
Moreover, the reduced vector is I-representable on $\Z$.
If $\sigma$ also represents an extension $(m,t)$,
then $\sigma_a$ represents $M_a(m,t)$.

\item Suppose, in addition, that $n=2r-1$ and that a
minimizing polynomial factors as
$p_{2r}^{(m)}=B_aR$, with $R\in\cE_{2r-2}$.
Then $M_a(m,\tau_{\Z}(m))$ is the extremal $\Z-$completion
of $M_a(m)$, and
\begin{equation}
\label{eq:Z-reduced-extremal-support}
\supp\mu_{\Z}(M_a(m))
=\supp\mu_{\Z}(m)\setminus\{a,a+1\}.
\end{equation}
In particular, if both $a$ and $a+1$ belong to
$\supp\mu_{\Z}(m)$, the original support is recovered
by adjoining these two points to the reduced support.
\end{enumerate}
\end{lemma}

\begin{proof}
Let $a\in\Z$ and $\sigma$ be any $\Z-$representing
probability measure for $m$. Since $B_a\ge0$ on $\Z$ and
$L_m(B_a)>0$, the measure in
\eqref{eq:Z-weighted-measure} is a probability measure.
By construction,
\[
\int_{\Z}x^\ell\,d\sigma_a(x)
=\frac{L_m(x^\ell B_a)}{L_m(B_a)}
=M_{a,\ell},
\qquad 0\le\ell\le n-2.
\]
The same identity holds up to order $n-1$ whenever
$\sigma$ also represents $(m,t)$.

We now prove that $M_a(m)$ is I-representable on $\Z$.
Recall that $d_n=2\lfloor n/2\rfloor$, so that
$d_{n-2}=d_n-2$, and let $q\in\cE_{d_n-2}$.
Since $\sigma_a$ represents $M_a(m)$, we have
$L_{M_a(m)}(q)\ge0$. If equality held, then
\[
0=L_{M_a(m)}(q)
=\frac{1}{L_m(B_a)}
\int_{\Z}B_a(x)q(x)\,d\sigma(x).
\]

Since $B_aq$ is non-negative on $\mathbb{Z}$, this would imply
$\operatorname{supp}\sigma\subseteq Z(q)\cup\{a,a+1\}$.
Applying \thref{lem:Z-compatible-completion}\textup{(ii)} for $r=d_n/2$, $R=q\in\mathcal{E}_{2r-2}$, and $S:=\operatorname{supp}\sigma\setminus\{a,a+1\}$,
we obtain $P\in\mathcal{E}_{d_n}$ with
$\operatorname{supp}\sigma\subseteq S\cup\{a,a+1\}\subseteq Z(P)$.
Hence, $L_m(P)=\int_{\mathbb{Z}}P\,d\sigma=0$, which contradicts the
I-representability of $m$. This proves (i).

For \textup{(ii)}, consider $\sigma=\mu_{\Z}(m)$ and
$t=\tau_{\Z}(m)$. By \textup{(i)}, the probability
measure $\sigma_a$ represents $M_a(m,t)$, whose
odd-order prefix $M_a(m)$ is I-representable.
Since $p_{2r}^{(m)}=B_aR$ is a minimizing polynomial of order $2r$, we have
\[
L_{M_a(m,t)}(R)
=\frac{L_{(m,t)}(B_aR)}{L_m(B_a)}=0.
\]
Thus, $M_a(m,t)$ is B-representable, since
$R\in\cE_{2r-2}$. By
\thref{prop:Z-inductive}\textup{(ii)} it is,
therefore, the extremal $\Z-$completion of
$M_a(m)$, and its unique representing measure
is $\mu_{\Z}(M_a(m))=\sigma_a$.

Finally, $B_a$ vanishes precisely at $a$ and
$a+1$ and is strictly positive at every other
integer. The definition of $\sigma_a$, therefore,
gives
\[
\supp\sigma_a
=\supp\sigma\setminus\{a,a+1\},
\]
which proves \eqref{eq:Z-reduced-extremal-support}
and the last assertion.
\end{proof}

In part \textup{(ii)}, the condition $R\in\cE_{2r-2}$ matters:
removing two consecutive zeros belonging to different lattice blocks
need not leave a polynomial in this family. The lattice block supplied by
\thref{lem:Z-forced-block} satisfies the required factorization.

\section{Algorithms}
\label{sec:decision-Z}

The preceding results yield two finite procedures: \textnormal{\textsc{MinZ}}
computes extremal $\Z-$completions and their representing measures,
while \textnormal{\textsc{DecideZ}} tests $\Z-$representability.

\subsection{The subroutine \textnormal{\textsc{MinZ}}}
The first procedure computes the extremal completion by recursively
reducing the order of the input moment vector.

\noindent\textbf{Input.} An I-representable vector $m=m^{(2r-1)}$.

\noindent\textbf{Output.} A minimizing polynomial $p\in\cE_{2r}$,
the extremal $\Z-$representing measure $\mu_{\Z}(m)$, and its support $S$.
\begin{enumerate}[label=\textbf{\arabic*.},leftmargin=0.8cm,itemsep=3pt]
\item \emph{Base case.} If $r=1$, set $a:=\lfloor m_1\rfloor$, return $p=B_a$
and the measure in \eqref{eq:Z-order-two}
with $S$ consisting only of the points with positive weight.
\item \emph{Integer support case.}
If $r\ge2$, compute the extremal Hamburger measure
$\mu_{\R}(m)$ as described in
\thref{prop:Z-hamburger-completion},
using \eqref{eq:Z-hamburger-system}. If its support is contained in $\Z$, use
the greedy construction in \thref{lem:Z-compatible-completion}\textup{(i)}
to choose a minimizing polynomial $p\in\cE_{2r}$. It is sufficient to choose one whose zero set contains $S:=\supp\mu_\mathbb{R}(m)$.
Return $(p,\mu_\mathbb{R}(m),S)$.
\item \emph{Recursive case.} Otherwise, initialize $\mathcal C:=\varnothing$. For each
$a\in\mathcal A_H(m)$, defined in \eqref{eq:Z-admissible-blocks}, compute $M_a(m)$ and call
$\textnormal{\textsc{MinZ}}(M_a(m))$, obtaining $(p_a,\mu_a,S_a)$.
If $S_a\cap\{a,a+1\}\ne\varnothing$, discard this branch. If $S_a\cap\{a,a+1\} = \varnothing$, use the greedy construction in
\thref{lem:Z-compatible-completion}\textup{(ii)} to choose
$q_a\in\cE_{2r}$ with
$S_a\cup\{a,a+1\}\subseteq Z(q_a)$, and add $q_a$ to $\mathcal C$.
\item \emph{Selection and reconstruction.} Choose $p\in\mathcal C$ which minimizes $\Phi_m$.
Write $Z(p)=\{z_1<\cdots<z_{2r}\}$ and solve
\begin{equation}
\label{eq:Z-algorithm-vandermonde}
\sum_{j=1}^{2r}\lambda_jz_j^\ell=m_\ell,
\qquad \ell=0,\ldots,2r-1.
\end{equation}
Return $\bigl(p,\sum_j\lambda_j\delta_{z_j},
\{z_j:\lambda_j>0\}\bigr)$.
\end{enumerate}

\noindent\emph{Correctness and termination.}
The correctness of \textnormal{\textsc{MinZ}} follows by induction on $r$.
The case $r=1$ is given by \eqref{eq:Z-order-two}, while integer Hamburger
nodes are covered by \thref{prop:Z-hamburger-completion}. Otherwise, every
recursive input is I-representable by
\thref{lem:Z-algorithm-block-reduction}\textup{(i)}. By
\thref{lem:Z-forced-block} and
\thref{lem:Z-algorithm-block-reduction}\textup{(ii)}, at least one branch
survives and yields a candidate $q_a$ vanishing on
$\supp\mu_{\Z}(m)$. Thus, $\mathcal C$ is non-empty and the selected
$p$ satisfies
\[
0\le L_{(m,\tau_{\Z}(m))}(p)
\le L_{(m,\tau_{\Z}(m))}(q_a)=0.
\]
Hence, $p$ is a minimizing polynomial, and
\eqref{eq:Z-algorithm-vandermonde} recovers the extremal
$\Z-$representing measure with non-negative weights. Termination follows
because each recursive call reduces~$r$.

\subsection{The decision procedure \textnormal{\textsc{DecideZ}}}
The second procedure processes the moments in increasing order,
using \textnormal{\textsc{MinZ}} at each even-order interior extension.

\noindent\textbf{Input.} A vector $m^{(n)}\in\R^n$.

Initialize the status to $\operatorname{I}$ and
process $k=2,\ldots,n$ in increasing order.
\begin{enumerate}[label=\textbf{\arabic*.},leftmargin=0.8cm,itemsep=3pt]
\item \emph{Boundary status.} If the current status is $\operatorname{B}$, with stored measure
$\mu$, compare $m_k$ with $\int x^k\,d\mu$. If they differ, return
$\operatorname{REJECT}$; otherwise keep the status $\operatorname{B}$ and the measure,
and continue with the next value of $k$.
\item \emph{Interior status.} Otherwise the status is $\operatorname{I}$. If $k$ is odd,
continue with the next value of $k$. If $k=2r$, call $\textnormal{\textsc{MinZ}}(m^{(2r-1)})$ and let $(p,\mu_*,S_*)$
be its output. Set $\Delta:=L_{m^{(2r)}}(p)$.
For $\Delta<0$, return $\operatorname{REJECT}$; for $\Delta=0$, set the
status to $\operatorname{B}$ and store $\mu_*$; for $\Delta>0$, retain
$\operatorname{I}$.
\end{enumerate}
After processing all moments, return $\operatorname{I}$ or
$(\operatorname{B},\mu)$ according to the final status. When $n=1$, the
loop is empty and the output is $\operatorname{I}$.

\noindent\emph{Correctness and termination.}
By \thref{prop:Z-inductive}, \textnormal{\textsc{DecideZ}} correctly
distinguishes non-representable, I-representable, and B-representable
vectors: part \textup{(iii)} handles boundary prefixes, part
\textup{(i)} handles odd-order extensions of interior prefixes, and
part \textup{(ii)}, together with the correctness of
\textnormal{\textsc{MinZ}}, handles even-order extensions.
Since both the outer loop and the recursive subroutine terminate,
this gives a finite decision algorithm in exact arithmetic.
In the B-representable case, the stored measure is the unique
$\Z-$representing probability measure.

\begin{remark}[Complexity of the recursion]
Let $T(r)$ denote the maximum number of calls to
\textnormal{\textsc{MinZ}} on an input of order $2r-1$,
including the initial call. Then $T(1)=1$ and
\[
T(r)\le 1+rT(r-1),
\]
which yields $T(r)=\mathcal{O}(r!)$.
\end{remark}

\begin{example}[A fourth-order computation]
\thlabel{ex:Z-algorithm-order-four}
Let us consider $m^{(4)}=(\frac12,2,-6,t)$.
Since $B_0$ is a minimizing polynomial of order two
for $m_1=\frac12$, the second-order test gives
\[
L_{(\frac12,2)}(B_0)=\frac32>0.
\]
Hence, $(\frac12,2)$ is I-representable on $\Z$,
and so is $(\frac12,2,-6)$ by
\thref{prop:Z-inductive}\textup{(i)}.

For $m=(\frac12,2,-6)$, the Hamburger system gives
\[
A_1=
\begin{pmatrix}
1&\frac12\\
\frac12&2
\end{pmatrix},
\qquad
h_2=
\begin{pmatrix}
2\\-6
\end{pmatrix}, \qquad q_2(x)=x^2+4x-4, 
\qquad
\tau_{\R}(m)=32.
\]
The nodes
\[
y_1=-2-2\sqrt2\in(-5,-4),
\qquad
y_2=-2+2\sqrt2\in(0,1)
\]
yield $\mathcal A_H(m)=\{-5,0\}$.
The two branches are
\[
\begin{array}{c|c|c|c}
a&L_m(B_a)&M_{a,1}
&\supp\mu_{\Z}(M_a(m))\\ \hline
-5&\frac{53}{2}&\frac{44}{53}&\{0,1\}\\
 0&\frac32&-\frac{16}{3}&\{-6,-5\}
\end{array}
\]
and both survive the support-disjointness test.
Indeed,
\[
\{0,1\}\cap\{-5,-4\}=\varnothing,
\qquad
\{-6,-5\}\cap\{0,1\}=\varnothing.
\]
The recursive base cases return $B_0$ and
$B_{-6}$, respectively. Consequently,
the two branches yield the candidates
\[
\begin{aligned}
p_{-5}
&=B_{-5}B_0=(x+5)(x+4)x(x-1)=x^4+8x^3+11x^2-20x,
\\[2mm]
p_0
&=B_0B_{-6}=x(x-1)(x+6)(x+5)=x^4+10x^3+19x^2-30x.
\end{aligned}
\]
Their values under the functional of \thref{def:Z-minimizing-polynomial} are given by
\[
\begin{aligned}
\Phi_m(p_{-5})
&=8(-6)+11(2)-20\left(\frac12\right)
=-36,\\
\Phi_m(p_0)
&=10(-6)+19(2)-30\left(\frac12\right)
=-37,
\end{aligned}
\]
and since $-37<-36$, the second candidate is
selected. Thus, $p_4^{(m)}
=B_0B_{-6}$
is a minimizing polynomial of order $4$.
In particular,
$L_{m^{(4)}}(p_4^{(m)})=t-37$, and
\[
\textnormal{\textsc{DecideZ}}
\left(\frac12,2,-6,t\right)=
\begin{cases}
\operatorname{REJECT},&t<37,\\
(\operatorname{B},\mu),&t=37,\\
\operatorname{I},&t>37,
\end{cases}
\]
where
\[
\mu=
\frac1{84}\delta_{-6}
+\frac1{30}\delta_{-5}
+\frac{13}{60}\delta_0
+\frac{31}{42}\delta_1.
\]
In particular, the lattice constraint raises the
smallest admissible fourth moment from $32$
to $37$. In the boundary case, further moments are forced by
the unique $\Z-$representing probability measure
$\mu$; for example, a $\Z-$representable extension
to order six must have $m_5=-196$ and $m_6=1077$.
\end{example}

\clearpage
\appendix
\numberwithin{equation}{section}
\section{Core technical results}
\label{app:core-technical-results}

We now establish the finite-reduction results used in the proof of
\thref{thm:Z-main}. The argument adapts the scale-reduction method of
\cite[Theorem~3.2]{InfusinoKunaLebowitzSpeer2017} to the full integer lattice $\Z$.
The key step is to remove the dependence on \(N\), which is achieved in the next results.
Recall that $\rho(p)$ denotes the largest absolute value
of a root of $p$, with the convention $\rho(1)=0$
(see \eqref{eq:rho(P)-def}).

\subsection{Scale classes and canonical displacements}

Given a polynomial $p\in\cE_{2r}$ we order its lattice blocks lexicographically by
$(\sigma(a),a)$, where \mbox{$\sigma(a):=\max\{|a|,|a+1|\}$}, and we write
\[
p=\prod_{h=1}^r B_{a_h},\qquad A_h(p):=\sigma(a_h).
\]
Thus, $A_1(p)\le\cdots\le A_r(p)=\rho(p)$ and ordering by $a$ resolves
ties between lattice blocks of equal scale.

\begin{definition}[Scale classes]
\thlabel{def:app-scale-classes}
Let $\mathcal J=(J_1,\ldots,J_r)$ be a strictly increasing tuple of
positive integers. For $0\le\ell\le r$, we define
\[
\cE_{2r,\ell,\mathcal J}:=
\left\{p\in\cE_{2r}:
 A_\ell(p)\le J_\ell\ \text{if }\ell>0,\quad
 A_{\ell+1}(p)>J_{\ell+1}\ \text{if }\ell<r\right\}.
\]
\end{definition}

The classes $\cE_{2r,\ell,\mathcal J}$
cover $\cE_{2r}$, although they need
not be disjoint. Indeed, for any $p\in\cE_{2r}$, it suffices to take the
largest index $\ell$ such that $A_\ell(p)\le J_\ell$,
or $\ell=0$ if no such index exists.
Moreover, the class $\cE_{2r,r,\mathcal J}$ is finite,
since it consists of polynomials whose roots belong
to $K_{J_r}$.

For $p\in\cE_{2r,\ell,\mathcal J}$ with $\ell<r$, let
$q:=\prod_{h=1}^{\ell}B_{a_h}$ with the convention $q=1$ when $\ell=0$.
We call these the small lattice blocks; their roots lie in $K_{J_\ell}$ when
$\ell>0$, while the remaining large lattice blocks have scale greater than
$J_{\ell+1}$. Let $r_-$ and
$r_+$ be their numbers on the negative and positive sides of the integer lattice $\Z$, respectively.
We define the corresponding canonical polynomial as
\begin{equation}
\label{eq:app-canonical-template}
p^*:=q
\prod_{h=1}^{r_-}B_{-J_{\ell+1}-2h+1}
\prod_{h=1}^{r_+}B_{J_{\ell+1}+2h-2},
\qquad r_-+r_+=r-\ell.
\end{equation}
On each side of the origin, the added lattice blocks are placed as close to the origin as possible, subject to having scale greater than \(J_{\ell+1}\).

\begin{definition}[Canonical scale classes]
\thlabel{def:app-canonical-scale-classes}
For $0\le\ell<r$, let $\cE^*_{2r,\ell,\mathcal J}$
be the set of all polynomials of the form
\eqref{eq:app-canonical-template}, where
$q\in\cE_{2\ell}$ and $\rho(q)\le J_\ell$
if $\ell>0$, while $q=1$ if $\ell=0$.
We also set $\cE^*_{2r,r,\mathcal J}:=\cE_{2r,r,\mathcal J}.$
\end{definition}

Every canonical scale class $\cE^*_{2r,\ell,\mathcal J}$
is finite. Indeed, there are finitely many admissible
small factors $q$ and, for each of them, only
$r-\ell+1$ possible choices for $(r_-,r_+)$.

\begin{definition}[Inward displacement at a threshold]
\thlabel{def:app-inward-displacement}
Let $q\in\cE_{2\ell}$ with $0\le\ell<r$, and let $M>\rho(q)$ be an integer.
An \emph{inward displacement at threshold $M$ fixing $q$} is a move
between polynomials in $\cE_{2r}$ that leaves the factor $q$ unchanged
and replaces one of the remaining lattice blocks according to
\[
B_a\longmapsto B_{a-1}\quad\text{if $a\ge0$},\qquad
B_a\longmapsto B_{a+1}\quad\text{if $a\le-1$},
\]
with all lattice blocks outside $q$ having scale greater than $M$
both before and after the move.
\end{definition}

\begin{lemma}[Canonical displacement]
\thlabel{lem:app-canonical-displacement}
For every $p\in\cE_{2r,\ell,\mathcal J}$ with $\ell<r$, take $q$ to be its
small factor and $M=J_{\ell+1}$. Then $p$ can be transformed into its
canonical polynomial $p^*$ by finitely many inward displacements at
threshold $M$ fixing $q$. Throughout the construction, the polynomial
remains in $\cE_{2r,\ell,\mathcal J}$ and its root radius does not
increase. In particular, $\rho(p^*)\le\rho(p)$.
\end{lemma}

\begin{proof}
We first write $p$ as in \eqref{eq:E-lattice-blocks}. On either side of the lattice, we list the large lattice blocks from the origin outwards. Their
endpoints satisfy
\[
a_h^+\ge J_{\ell+1}+2h-2,\qquad
a_h^-\le-J_{\ell+1}-2h+1,
\]
because each lattice block has scale greater than $J_{\ell+1}$ and distinct lattice blocks
have disjoint roots. We move the innermost lattice block to its canonical position,
then the next lattice block, and so on, treating the two sides separately.
There is no collision with the small roots, which lie strictly between the
two groups of large lattice blocks. Nor is there a collision with another large
lattice block: lattice blocks already moved occupy their inner canonical positions, whereas
the lattice blocks not yet moved lie farther out. Every move decreases the distance
to the prescribed endpoint by one, preserves the scale bounds, and moves
roots towards the origin. The procedure therefore terminates at $p^*$.
\end{proof}
\begin{example}
    Let $r=4$, $\ell=1$, $\mathcal J=(2,5,9,14)$, and consider $p=B_{-1}B_{-9}B_8B_{12}.$
The corresponding lattice-block scales are $1,9,9,13$, respectively.
Since $A_1(p) = 1 \leq 2 = J_1$ and $A_2(p) = 9 > 5 = J_2$, we have $p\in\cE_{8,1,\mathcal J}$. Moreover, its small factor is $q=B_{-1}$,
while $(r_-,r_+)=(1,2)$. The corresponding canonical polynomial is, therefore, $p^*=B_{-1}B_{-6}B_5B_7.$
We first move the innermost positive-side lattice block to its
canonical position,
\[
B_8\longmapsto B_7\longmapsto B_6\longmapsto B_5,
\]
and only then move the outermost positive-side lattice block:
\[
B_{12}\longmapsto B_{11}\longmapsto B_{10}
\longmapsto B_9\longmapsto B_8\longmapsto B_7.
\]
Finally, on the negative side, we perform
\[
B_{-9}\longmapsto B_{-8}\longmapsto B_{-7}\longmapsto B_{-6}.
\]

The order matters: if we moved $B_{12}$ first while keeping $B_8$
fixed, the step $B_{10}\longmapsto B_9$ would produce a factor
$B_8B_9$ with a repeated root at $9$, because
\[
B_8(x)=(x-8)(x-9),\qquad B_9(x)=(x-9)(x-10).
\]
The resulting polynomial would therefore not belong to $\cE_8$.
Thus, by moving the innermost lattice block first, we keep all roots distinct, which guarantees that every intermediate polynomial belongs to $\cE_{8,1,\mathcal J}$,
and the root radius never increases. Figure \ref{fig:canonical-displacement} illustrates the initial and the canonical configurations respectively, making clear how the lattice blocks are displaced inward while preserving the fixed factor and avoiding overlapping roots.
\begin{figure}[H]
\centering
\begin{tikzpicture}[x=0.51cm,y=0.95cm,>=Stealth]
  \def\ytop{2.4}
  \def\ybottom{0}
  \fill[blue!6] (-2,-0.52) rectangle (2,3.51);
  \foreach \x in {-2,2}
    \draw[blue!65!black,densely dotted] (\x,-0.48) -- (\x,3.57);
  \foreach \x in {-5,5}
    \draw[orange!80!black,dashed] (\x,-0.48) -- (\x,3.57);
  \node[above,font=\scriptsize,text=blue!65!black]   at (-2,3.65) {$-J_1=-2$};
  \node[above,font=\scriptsize,text=blue!65!black]   at ( 2,3.65) {$J_1=2$};
  \node[above,font=\scriptsize,text=orange!80!black] at (-5,3.65) {$-J_2=-5$};
  \node[above,font=\scriptsize,text=orange!80!black] at ( 5,3.65) {$J_2=5$};

  \foreach \y in {0,2.4}{
    \draw[-{Stealth[length=2.2mm]},thin] (-11,\y) -- (14,\y);
    \foreach \x in {-10,-8,...,12}{
      \draw (\x,{\y-0.075}) -- (\x,{\y+0.075});
      \node[below=3pt,font=\scriptsize] at (\x,{\y-0.07}) {$\x$};
    }
  }
  \node[left] at (-11.05,\ytop) {$Z(p)$};
  \node[left] at (-11.05,\ybottom) {$Z(p^*)$};

  \newcommand{\pictureblock}[5]{%
    \draw[line width=1.3pt,#3] (#1,{#2+0.22}) -- ({#1+1},{#2+0.22});
    \fill[#3] (#1,{#2+0.22}) circle (1.65pt);
    \fill[#3] ({#1+1},{#2+0.22}) circle (1.65pt);
    \node[above=3pt,font=\scriptsize,align=center,text=#3]
      at ({#1+0.5},{#2+0.22}) {$#4$\\[-2pt]$\sigma=#5$};
  }
  \pictureblock{-9}{2.4}{black}{B_{-9}}{9}
  \pictureblock{-1}{2.4}{blue!75!black}{B_{-1}}{1}
  \pictureblock{8}{2.4}{black}{B_8}{9}
  \pictureblock{12}{2.4}{black}{B_{12}}{13}

  \pictureblock{-6}{0}{black}{B_{-6}}{6}
  \pictureblock{-1}{0}{blue!75!black}{B_{-1}}{1}
  \pictureblock{5}{0}{black}{B_5}{6}
  \pictureblock{7}{0}{black}{B_7}{8}

  \draw[-{Stealth[length=2mm]},line width=0.7pt] (3.1,1.72) -- (3.1,0.70);
\end{tikzpicture}
\caption{Canonical displacement for $\ell=1$ and $\mathcal J=(2,5,9,14)$.
The dotted and dashed lines mark the thresholds $\pm J_1$ and $\pm J_2$, respectively.
The labels show the block scales $\sigma(a)=\max\{|a|,|a+1|\}$.
The small factor $q=B_{-1}$ (blue) remains fixed, while the remaining blocks
are displaced inward to obtain $p^*=B_{-1}B_{-6}B_5B_7$.}
\label{fig:canonical-displacement}
\end{figure}
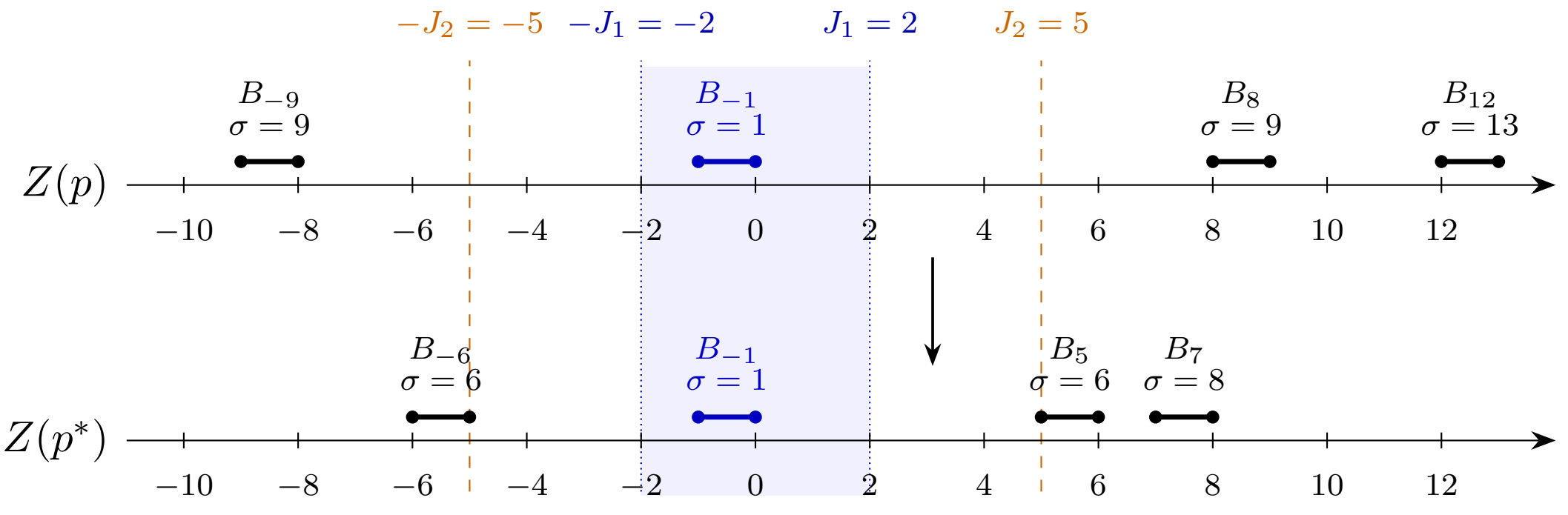
\end{example}

\subsection{Coefficient estimates and propagation of zeros}

The following estimates provide the uniform bounds for sufficiently distant roots required by the domination argument used in Subsection \ref{subsec:uniform-scale-reduction}.

\begin{definition}[Elementary symmetric polynomials]
Let $y_1,\ldots,y_s\in\mathbb{R}$. For $0\le j\le s$, the
$j$-th elementary symmetric polynomial in $y_1,\ldots,y_s$ is defined by
\[
e_j(y_1,\ldots,y_s)
:=
\sum_{1\le i_1<\cdots<i_j\le s}
y_{i_1}\cdots y_{i_j}
\]
with the convention $e_0(y_1,\ldots,y_s):=1.$
\end{definition}

\begin{lemma}[Coefficient estimates]
\thlabel{lem:app-coefficient-estimates}
Let $G(x)=\sum_{k=0}^d g_kx^k$ be a polynomial of degree $d$ with $G(0)=g_0>0$
and roots $\tau_1,\ldots,\tau_d$ with $|\tau_i|\ge M>0$. Then
\begin{equation}
\label{eq:app-coefficient-bound}
\frac{|g_k|}{g_0}
\le\binom{d}{k}M^{-k},
\qquad 1\le k\le d.
\end{equation}
\end{lemma}

\begin{proof}
Since none of the roots $\tau_1,\ldots,\tau_d$ is zero,
by Vi\`ete's formulas (see e.g. \cite{Lang2002}) we have
\[
g_k=(-1)^{d-k}g_d\,e_{d-k}(\tau_1,\ldots,\tau_d),
\qquad 0\le k\le d.
\]
In particular $g_0=(-1)^d g_d\prod_{i=1}^d\tau_i\ne0$, and consequently we get
\[
\begin{aligned}
\frac{g_k}{g_0}
&=(-1)^k
\frac{e_{d-k}(\tau_1,\ldots,\tau_d)}
{\prod_{i=1}^d\tau_i}=(-1)^k
\sum_{\substack{I\subseteq\{1,\ldots,d\}\\|I|=d-k}}
\frac{\prod_{i\in I}\tau_i}
{\prod_{i=1}^d\tau_i}=(-1)^k
\sum_{\substack{I\subseteq\{1,\ldots,d\}\\|I|=d-k}}
\prod_{j\notin I}\frac1{\tau_j}\\
&=(-1)^k
\sum_{\substack{J\subseteq\{1,\ldots,d\}\\|J|=k}}
\prod_{j\in J}\frac1{\tau_j}=(-1)^k
e_k\left(\frac1{\tau_1},\ldots,\frac1{\tau_d}\right),
\end{aligned}
\]
where the fourth equality follows by setting
$J=\{1,\ldots,d\}\setminus I$.
There are exactly $\binom{d}{k}$ terms in the last sum, and each
has absolute value at most $M^{-k}$. Hence, by the triangle inequality, we get
\[
\frac{|g_k|}{g_0}
\le\sum_{\substack{J\subseteq\{1,\ldots,d\}\\|J|=k}}
\prod_{j\in J}\frac1{|\tau_j|}
\le\binom{d}{k}M^{-k},
\]
which proves \eqref{eq:app-coefficient-bound}.
\end{proof}

\begin{lemma}[Propagation of zeros]
\thlabel{lem:local-scale-bounds}
Let $m\in\R^{2n}$ such that $L_m(p)\ge0$ for every
$p\in\cE_{2n}$. Suppose that $q\in\cE_{2\ell}$ for some
$0\le\ell<n$ and $L_m(q)=0$, and set $s:=2(n-\ell)$. Then:
\begin{enumerate}[label=\textup{(\alph*)},leftmargin=1cm]
\item $L_m(x^kq)=0$ for every $0\le k<s$, and $L_m(x^sq)\ge0$;
\item if, in addition,
\begin{equation}
\label{eq:app-zero-compatibility}
L_m(H)=0\implies L_m(x^2H)=0
\quad\text{for every}\quad H\in\cE_{2n-2},
\end{equation}
then also $L_m(x^sq)=0$.
\end{enumerate}
\end{lemma}
\begin{proof}
By repeatedly applying \thref{lem:positivityZ}, we obtain
$L_m(p)\ge0$ for every $p\in\cE_{2j}$, $0\le j\le n$.
We first establish an auxiliary observation: if
$H\in\cE_{2j}$, with $0\le j\le n-2$, satisfies
$L_m(H)=0$, then
\[
L_m(HB_a)=0
\]
for every lattice block $B_a$ whose roots are disjoint
from those of $H$. Indeed, by \thref{lem:positivityZ},
\[
L_m(xH)=0,\qquad L_m(x^2H)\ge0.
\]
For all sufficiently large integers $M$, the roots of
$B_{-M-1}B_M$ are disjoint from those of $H$ and, therefore,
$HB_{-M-1}B_M\in\cE_{2j+4}$ and its Riesz value is
non-negative. Since
\[
B_{-M-1}(x)B_M(x)
=x^4-(2M^2+2M+1)x^2+M^2(M+1)^2,
\]
and $L_m(H)=0$, we obtain
\[
0\le L_m(x^4H)-(2M^2+2M+1)L_m(x^2H).
\]
The first term does not depend on $M$ and, thus, letting $M$ tend
to infinity gives $L_m(x^2H)\le0$. Consequently,
$L_m(x^2H)=0$ and for every lattice block $B_a$ whose roots
are disjoint from those of $H$, we then have
\[
\begin{aligned}
L_m(HB_a)
&=L_m(x^2H)-(2a+1)L_m(xH)+a(a+1)L_m(H)=0.
\end{aligned}
\]
We now apply this observation to $q$. Choose
$n-\ell-1$ lattice blocks
$B_{a_1},\ldots,B_{a_{n-\ell-1}}$ with pairwise
disjoint roots, also disjoint from those of $q$.
We define
\[
p_0:=1,\qquad
p_j:=\prod_{i=1}^j B_{a_i},
\quad 1\le j\le n-\ell-1,
\]
and we claim that
\[
L_m(qp_j)=0,\qquad 0\le j\le n-\ell-1.
\]
The case $j=0$ is our assumption. If the claim holds
for some $j\le n-\ell-2$, then
$qp_j\in\cE_{2\ell+2j}$ has degree at most $2n-4$.
The preceding observation therefore gives
\[
L_m(qp_{j+1})
=L_m(qp_jB_{a_{j+1}})=0,
\]
which concludes the induction step.
Moreover, since $qp_j\in\cE_{2\ell+2j}$ and
$2\ell+2j\le2n-2$, by \thref{lem:positivityZ} we have
\[
L_m(xqp_j)=0,\qquad 0\le j\le n-\ell-1.
\]
Consequently, $L_m(qh)=0$ whenever $h$ belongs
to the family
\[
p_0,\ xp_0,\ p_1,\ xp_1,\ \ldots,\
p_{n-\ell-1},\ xp_{n-\ell-1}.
\]
These polynomials are monic of respective degrees
$0,1,\ldots,2(n-\ell)-1$ and, therefore, form a basis
of $\R[x]_{\le 2(n-\ell)-1}$. By linearity,
\[
L_m(x^kq)=0,\qquad 0\le k<2(n-\ell).
\]
Finally, consider the polynomial $H:=qp_{n-\ell-1}\in\cE_{2n-2}$.
Since $L_m(H)=0$, \thref{lem:positivityZ} gives
$L_m(x^2H)\ge0$. Moreover, since $p_{n-\ell-1}$
is monic of degree $2(n-\ell)-2$, we can write
\[
x^2p_{n-\ell-1}(x)
=x^{2(n-\ell)}
+\sum_{k=2}^{2(n-\ell)-1}c_kx^k,
\]
where the sum is empty if $n-\ell=1$.
Therefore,
\[
\begin{aligned}
L_m(x^2H)
&=L_m(x^{2(n-\ell)}q)
  +\sum_{k=2}^{2(n-\ell)-1}c_kL_m(x^kq)=L_m(x^{2(n-\ell)}q),
\end{aligned}
\]
where the last equality follows from the vanishing
identities established above. Hence, $L_m(x^{2(n-\ell)}q)=L_m(x^2H)\ge0,$
which completes \textup{(a)}.
If \eqref{eq:app-zero-compatibility} also holds,
then $L_m(x^2H)=0$, proving \textup{(b)}.
\end{proof}

\subsection{Uniform scale reduction}
\label{subsec:uniform-scale-reduction}

The polynomial families considered in \thref{thm:Z-main} are infinite.
We show that, under the hypotheses of \thref{lem:reduction}, inward displacements
reduce the relevant Riesz evaluations to those of a finite family of
canonical polynomials.

\begin{lemma}[One-step difference]
\thlabel{lem:app-one-step-displacement}
Let $p=qB_aR\in\cE_{2r}$, where $0\le\ell<r$,
$q\in\cE_{2\ell}$ is a product of $\ell$ disjoint lattice blocks
of $p$, and $R$ is the product of the remaining $r-\ell-1$
blocks (with $R=1$ if $\ell=r-1$).
Let $M>\rho(q)$ be an integer and suppose that $p'$ is obtained
from $p$ by an inward displacement at threshold $M$ fixing $q$,
as in \thref{def:app-inward-displacement}, with $B_a$ the displaced
block. In particular, all blocks of $B_aR$ have scale greater
than $M$ before the displacement, and the displaced block still
has scale greater than $M$ afterwards. Then
\begin{equation}
\label{eq:app-one-step-difference}
p-p'=qD,\qquad
D=\begin{cases}
-2(x-a)R,&\text{if }B_a\longmapsto B_{a-1},\\
2(x-a-1)R,&\text{if }B_a\longmapsto B_{a+1}.
\end{cases}
\end{equation}
Moreover, $\deg D=2(r-\ell)-1$, $D(0)>0$, and every root of $D$
has absolute value at least $M$.
\end{lemma}
\begin{proof}
The identities
\[
B_a-B_{a-1}=-2(x-a),\qquad B_a-B_{a+1}=2(x-a-1)
\]
give \eqref{eq:app-one-step-difference}. Since $R$ consists of
$r-\ell-1$ lattice blocks, we obtain $\deg D=2(r-\ell)-1$.
If $B_a\longmapsto B_{a-1}$, then $a\ge0$ and the new block
has scale $\sigma(a-1)>M>0$. As $M>\rho(q)\ge0$ is an integer,
we have $M\ge1$; hence $\sigma(a-1)=a>M$.
The linear factor $-2(x-a)$ is therefore positive at zero and
has its root outside $[-M,M]$.
If $B_a\longmapsto B_{a+1}$, then $a\le-1$ and
$\sigma(a+1)=-a-1>M > 0$. Thus, $a < -1$ and so $2(x-a-1)$ is positive at zero
and its root also lies outside $[-M,M]$.
Every block in $R$ has scale greater than $M$. Its two integer
roots consequently have absolute value at least $M$, and its
value at zero is positive. Hence, $R(0)>0$ and all its roots
have absolute value at least $M$. The same properties hold
for $D$, as claimed.
\end{proof}

\begin{lemma}[One-step estimate]
\thlabel{lem:app-one-step-estimate}
Under the assumptions of \thref{lem:app-one-step-displacement}, let us fix
\[
w\in\{1,\,N-x,\,N+x,\,N^2-x^2\},\qquad G:=wD,\qquad d:=\deg G.
\]
When $w\ne1$, take any integer $N\ge\rho(p)$. Then $G(0)>0$ and
\begin{equation}
\label{eq:app-one-step-lower-bound}
L_m\bigl(w(p-p')\bigr)
\ge G(0)\left(
L_m(q)-\sum_{k=1}^{d}\binom{d}{k}M^{-k}|L_m(x^kq)|
\right).
\end{equation}
If $L_m(x^kq)=0$ for $0\le k\le d$, then
$L_m\bigl(w(p-p')\bigr)=0$.
\end{lemma}

\begin{proof}
By \thref{lem:app-one-step-displacement}, $D(0)>0$ and every
root of $D$ has absolute value at least $M$.
If $w=1$, the same properties hold for $G=D$.
Otherwise, $w(0)>0$ and every root of $w$ has absolute value
$N\ge\rho(p)>M$. Thus, in all cases $G(0)>0$ and every root
of $G$ has absolute value at least $M$.
By \eqref{eq:app-coefficient-bound} we get
\begin{equation}
\label{eq:app-one-step-coefficients}
G(x)=\sum_{k=0}^{d}g_kx^k,\qquad g_0=G(0)>0,\qquad
\frac{|g_k|}{g_0}\le\binom{d}{k}M^{-k}\quad(1\le k\le d).
\end{equation}
Since $w(p-p')=qG$, by the linearity of $L_m$ we have
\[
\begin{aligned}
L_m\bigl(w(p-p')\bigr)
&=g_0L_m(q)+\sum_{k=1}^{d}g_kL_m(x^kq)\ge g_0L_m(q)-\sum_{k=1}^{d}|g_k|\,|L_m(x^kq)|\\
& = g_0\left(L_m(q) - \sum_{k=1}^{d}\frac{|g_k|}{g_0}\,|L_m(x^kq)|\right)
\\
&\ge G(0)\left(
L_m(q)-\sum_{k=1}^{d}\binom{d}{k}M^{-k}|L_m(x^kq)|
\right),
\end{aligned}
\]
where the last inequality comes from \eqref{eq:app-one-step-coefficients}. This proves \eqref{eq:app-one-step-lower-bound}, and the same expansion
gives zero when all the stated values vanish.
\end{proof}
For a fixed $m$, the expression in parentheses in
\eqref{eq:app-one-step-lower-bound}
depends only on $q$, $M$, and $d$.
In particular, it is independent of $N$ and of the positions
of all the lattice block factors outside $q$.
When $w=N\pm x$, it is also independent of the $\pm$ sign.

\begin{theorem}[Scale reduction]
\thlabel{lem:reduction}
Let $n\ge1$. The following assertions hold.
\begin{enumerate}[label=\textup{(\alph*)},leftmargin=1cm]
\item Let $m\in\R^{2n+1}$ satisfy $L_m(q)>0$ for every
$q\in\cE_{2j}$ with $0\le j\le n$, and set
$\widehat m:=(m_1,\ldots,m_{2n+1},0)$.
There exists a strictly increasing tuple of positive integers
$\mathcal J=(J_1,\ldots,J_{n+1})$ such that,
for $0\le\ell\le n$ and $p\in\cE_{2n+2,\ell,\mathcal J}$,
its canonical polynomial satisfies
\[
L_{\widehat m}(p^*)\le L_{\widehat m}(p).
\]
\item Let $m\in\R^{2n+1}$ satisfy $L_m(q)\ge0$ for every
$q\in\cE_{2n}$ and the condition \eqref{eq:app-zero-compatibility}.
There exists a strictly increasing tuple of positive integers
$\mathcal J=(J_1,\ldots,J_n)$ such that,
for $0\le\ell<n$ and $p\in\cE_{2n,\ell,\mathcal J}$,
\[
L_m((N\pm x)p^*)\le L_m((N\pm x)p)
\qquad\text{for every integer }N\ge\rho(p),
\]
with the same canonical polynomial for both signs.
\item Let $n\ge2$ and $m\in\R^{2n}$ satisfy $L_m(q)\ge0$ for every
$q\in\cE_{2n}$.
There exists a strictly increasing tuple of positive integers
$\mathcal J=(J_1,\ldots,J_{n-1})$ such that,
for $0\le\ell<n-1$ and $p\in\cE_{2n-2,\ell,\mathcal J}$,
\[
L_m((N^2-x^2)p^*)\le L_m((N^2-x^2)p)
\qquad\text{for every integer }N\ge\rho(p).
\]
\end{enumerate}
In all cases, $p^*$ belongs to the corresponding finite family
$\cE^*_{\deg p,\ell,\mathcal J}$ and satisfies $\rho(p^*)\le\rho(p)$.
\end{theorem}

\begin{proof}
We construct the tuple separately in each case. The number of lattice blocks
is $r=n+1$, $r=n$, or $r=n-1$, and the weight is $w=1$, $w=N\pm x$,
or $w=N^2-x^2$, in cases \textup{(a)}, \textup{(b)}, and
\textup{(c)}, respectively. For a small factor $q\in\cE_{2\ell}$,
the degree in \thref{lem:app-one-step-estimate} is given by
\begin{equation}
\label{eq:app-displacement-degree}
d=2(r-\ell)-1+\deg w
=
\begin{cases}
2n-2\ell+1,&\text{in \textup{(a)}},\\
2n-2\ell,&\text{in \textup{(b)}},\\
2n-2\ell-1,&\text{in \textup{(c)}}.
\end{cases}
\end{equation}
At each scale $0\le\ell<r$, every $q\in\cE_{2\ell}$ satisfies
one of the following alternatives:
\begin{equation}
\label{eq:app-local-alternative}
L_m(q)>0
\qquad\text{or}\qquad
L_m(x^kq)=0\quad\text{for all }0\le k\le d.
\end{equation}
In fact, in case \textup{(a)} we have $L_m(q)>0$ by assumption. In cases
\textup{(b)} and \textup{(c)}, by repeatedly applying
\thref{lem:positivityZ} we get $L_m(q)\ge0$.
If $L_m(q)=0$, we use \thref{lem:local-scale-bounds}
to obtain the required vanishing property, otherwise $L_m(q) > 0$. More precisely,
for case \textup{(b)}, part \textup{(b)} of that lemma
applied to the even truncation $(m_1,\ldots,m_{2n})$,
gives $L_m(x^kq)=0$ for all $0\le k\le 2n-2\ell$.
For case \textup{(c)}, part \textup{(a)} of the lemma gives
$L_m(x^kq)=0$ for all $0\le k\le 2n-2\ell-1$.
Thus, in both cases, the vanishing holds up to
the required degree.

Now set $J_0:=0$. We construct $J_1,\ldots,J_r$ by induction on the scale $\ell$.
The induction hypothesis at scale $\ell$ is that every inward displacement
at threshold $J_{\ell+1}$ fixing a factor $q\in\cE_{2\ell}$ with
$\rho(q)\le J_\ell$ satisfies
\begin{equation}
\label{eq:app-induction-property}
L_m\bigl(w(p-p')\bigr)\ge0,
\end{equation}
where in the weighted cases $w\neq 1$, this must hold for every integer $N\ge\rho(p)$.

\smallskip
\noindent\emph{Base case: $\ell=0$.}
The only small factor is $q=1$ and $L_m(1)=1$. With $d$ given by
\eqref{eq:app-displacement-degree} at $\ell=0$, we choose a positive
integer $J_1$ such that
$$
\sum_{k=1}^{d}\binom{d}{k}J_1^{-k}|L_m(x^k)|<L_m(1) = 1.
$$
This is possible because the sum tends to zero as $J_1\to\infty$.
The right-hand side of \eqref{eq:app-one-step-lower-bound} with
$q=1$ and $M=J_1$, is then strictly positive. Thus,
\eqref{eq:app-induction-property} holds at scale zero.

\smallskip
\noindent\emph{Inductive step: from $\ell-1$ to $\ell$.}
Let $1\le\ell<r$. Assume that
$0=J_0<J_1<\cdots<J_\ell$ have been already chosen and that
\eqref{eq:app-induction-property} holds at every scale
$0,\ldots,\ell-1$. We keep these thresholds fixed.
There are only finitely many $q\in\cE_{2\ell}$ with
$\rho(q)\le J_\ell$, since their roots lie in $K_{J_\ell}$.
For each such $q$, we choose an integer $M_q>J_\ell$ as follows.
If $L_m(q)>0$, we choose $M_q$ so large that
$$
\sum_{k=1}^{d}\binom{d}{k}M_q^{-k}|L_m(x^kq)|<L_m(q).
$$
Such a choice exists since the sum approaches zero as $M_q\to\infty$, and the same strict inequality holds with $M_q$ replaced by any $M\ge M_q$.
If $L_m(q)=0$, the second alternative in
\eqref{eq:app-local-alternative} applies, and we take
$M_q:=J_\ell+1$. In the first case, every threshold $M\ge M_q$
gives a strictly positive one-step difference; in the second,
every threshold $M\ge M_q$ gives zero. Both assertions follow
from \thref{lem:app-one-step-estimate}. We now define
$$
J_{\ell+1}:=
\max_{\substack{q\in\cE_{2\ell}\\\rho(q)\le J_\ell}}M_q.
$$
This maximum is finite, and its index set is non-empty
because $J_\ell\ge\ell$ and the polynomial $q=\prod_{j=0}^{\ell-1}B_{-\ell+2j}$
belongs to $\cE_{2\ell}$ and satisfies $\rho(q)=\ell$.
 Moreover, we have $J_{\ell+1}>J_\ell$ and
$J_{\ell+1}\ge M_q$ for every possible small factor $q$. Hence,
\eqref{eq:app-induction-property} holds at scale $\ell$.
Thus, the requirement now holds at all scales
$0,\ldots,\ell$, which completes the inductive step. We have therefore obtained a
strictly increasing tuple $\mathcal J=(J_1,\ldots,J_r)$.

\smallskip
\noindent\emph{Passage to the canonical polynomial.}
Let us fix $p\in\cE_{2r,\ell,\mathcal J}$ with $\ell<r$. By
\thref{lem:app-canonical-displacement}, there exists a finite sequence
\[
p=p^{(0)}\longmapsto p^{(1)}\longmapsto\cdots
\longmapsto p^{(T)}=p^*
\]
of inward displacements at threshold $J_{\ell+1}$, all fixing
the same small factor $q$. If $w\neq 1$, fix any integer
$N\ge\rho(p)$. Since each inward displacement moves a lattice block toward
the origin, we have
\[
\rho(p^{(T)})\le \rho(p^{(T-1)})\le\cdots
\le\rho(p^{(1)})\le\rho(p^{(0)})=\rho(p)\le N.
\]
Thus, all the polynomials in the above chain have their roots in
$[-N,N]$, and the same $N$ can be used throughout the
entire sequence of displacements. By repeated applications of
\eqref{eq:app-induction-property} and by linearity of $L_m$, we find
\[
L_m\bigl(w(p-p^*)\bigr)
=\sum_{h=0}^{T-1}
L_m\bigl(w(p^{(h)}-p^{(h+1)})\bigr)\ge0.
\]
In \textup{(a)}, $w = 1$ and the leading terms of $p$ and $p^*$ cancel because $p,p^*$ are monic. Therefore,
\[
L_{\widehat m}(p)-L_{\widehat m}(p^*)=L_m(p-p^*)\ge0.
\]
In \textup{(b)} and \textup{(c)}, by linearity we get
$L_m(wp)-L_m(wp^*)\ge0$, as required.

The choices of thresholds are independent of $N$ and of the sign $\pm$
in $N\pm x$. So is the canonical construction and, therefore, the
same $p^*$ works for all these choices. Its membership in the
stated finite family and the bound $\rho(p^*)\le\rho(p)$ follow
from \thref{lem:app-canonical-displacement} and
\eqref{eq:app-canonical-template}.
\end{proof}

\subsection{Minimizing polynomials and finite boundary families}

The preceding construction replaces each scale class by a finite family, which allows us to consider a finite reduction procedure.

\begin{theorem}[Finite reduction]
\thlabel{thm:thm3.2-gen}
Let $n\in\N$.
\begin{enumerate}[label=\textup{(\alph*)},leftmargin=1cm]
\item Let $m^{(2n+1)}=(m_1,\ldots,m_{2n+1})\in\R^{2n+1}$ and assume
$L_{m^{(2n+1)}}(q)>0$ for every $q\in\bigcup_{j=0}^{n}\cE_{2j}$.
Then there exists a polynomial
$p_{2n+2}^{(m^{(2n+1)})}\in\cE_{2n+2}$ such that, for every $t\in\R$,
\[
L_{(m_1,\ldots,m_{2n+1},t)}
 \bigl(p_{2n+2}^{(m^{(2n+1)})}\bigr)
=\min_{p\in\cE_{2n+2}}L_{(m_1,\ldots,m_{2n+1},t)}(p).
\]
Thus, a minimizing polynomial of order $2n+2$ associated with
$m^{(2n+1)}$ exists in the sense of
\thref{def:Z-minimizing-polynomial}.
\item Let $m\in\R^{2n+1}$ satisfy $L_m(q)\ge0$ for every
$q\in\cE_{2n}$ and \eqref{eq:app-zero-compatibility}.
There exists a finite set $\widetilde{\cE}_{2n}^{m}\subseteq\cE_{2n}$
such that, for every $p\in\cE_{2n}$, there is a
$\widetilde p\in\widetilde{\cE}_{2n}^{m}$ with
$\rho(\widetilde p)\le\rho(p)$ satisfying
\[
L_m((N+\varepsilon x)\widetilde p)
\le L_m((N+\varepsilon x)p)
\]
for every integer $N\ge\rho(p)$ and every
$\varepsilon\in\{-1,1\}$. In particular, the same representative
works for both signs and for all such $N$.
\item Let $m\in\R^{2n}$ satisfy $L_m(q)\ge0$ for every
$q\in\cE_{2n}$.
There exists a finite set
$\widetilde{\cE}_{2n-2}^{m}\subseteq\cE_{2n-2}$ such that,
for every $p\in\cE_{2n-2}$, there is a
$\widetilde p\in\widetilde{\cE}_{2n-2}^{m}$ with
$\rho(\widetilde p)\le\rho(p)$ satisfying
\[
L_m((N^2-x^2)\widetilde p)\le L_m((N^2-x^2)p)
\]
for every integer $N\ge\rho(p)$.
\end{enumerate}
\end{theorem}

\begin{proof}
For a tuple $\mathcal J=(J_1,\ldots,J_r)$, we consider the finite set
\[
\mathcal F_{2r}(\mathcal J):=
\bigcup_{\ell=0}^r\cE^*_{2r,\ell,\mathcal J}.
\]
Every $p\in\cE_{2r}$ belongs to a scale class
$\cE_{2r,\ell,\mathcal J}$.
If $\ell<r$, replace $p$ by its canonical polynomial $p^*$;
otherwise, retain $p$ itself. In either case, the representative
belongs to $\mathcal F_{2r}(\mathcal J)$ and its root radius is no
larger than that of $p$.

For \textup{(a)}, we apply \thref{lem:reduction}\textup{(a)} and take
$r=n+1$. Every $p\in\cE_{2n+2}$ has a representative in
$\mathcal F_{2n+2}(\mathcal J)$ with no larger
$L_{(m_1,\ldots,m_{2n+1},0)}$-value. Hence, the minimum is attained
on this finite set. Since all polynomials in $\cE_{2n+2}$ are monic,
replacing the last moment $0$ by $t$ adds the same constant $t$ to
every value, so the same polynomial minimizes for every $t\in\R$.

For \textup{(b)}, we use the tuple from
\thref{lem:reduction}\textup{(b)} and set $\widetilde{\cE}_{2n}^{m}:=\mathcal F_{2n}(\mathcal J).$
For each $p$, the representative described above depends only on
$p$ and the chosen tuple, not on $N$ or the sign.
By \thref{lem:reduction}\textup{(b)}, it does not increase the Riesz functional for any integer $N\ge\rho(p)$.
Thus, a single finite family and the same representative give
both inequalities simultaneously.

Finally, for \textup{(c)} with $n\ge2$, we apply
\thref{lem:reduction}\textup{(c)} and set
$\widetilde{\cE}_{2n-2}^{m}:=\mathcal F_{2n-2}(\mathcal J)$.
The same construction gives a representative that works for every
integer $N\ge\rho(p)$. When $n=1$, we can simply take
$\widetilde{\cE}_0^m:=\{1\}$.
\end{proof}

\section*{Acknowledgments}
M. Infusino is a member of GNAMPA group of INdAM and is indebted to the Fondazione di Sardegna for the support received within the project ENTROPIA (grant number F83C26000350007).

\section*{Declaration of generative AI and AI-assisted technologies in the manuscript preparation process}
During the preparation of this work, the authors used OpenAI GPT-6 and GPT-6.2 Astra and
Anthropic Claude Fable 5.2 as well as Opus 5.5 to improve the clarity and the readability of the manuscript. After using
these tools, the authors reviewed and edited the content as needed
and take full responsibility for the content of the published article.


\bibliographystyle{plain}
\bibliography{references}

@article{AlOshAlzaid1987,
  author    = {Al-Osh, M. A. and Alzaid, A. A.},
  title     = {First-order integer-valued autoregressive ({INAR}(1)) process},
  journal   = {Journal of Time Series Analysis},
  volume    = {8},
  number    = {3},
  pages     = {261--275},
  year      = {1987},
  doi       = {10.1111/j.1467-9892.1987.tb00438.x}
}

@article{CurtoFialkow1996,
  author    = {Curto, R. E. and Fialkow, L. A.},
  title     = {Solution of the truncated complex moment problem for flat data},
  journal   = {Memoirs of the American Mathematical Society},
  volume    = {119},
  number    = {568},
  pages     = {x+52},
  publisher = {American Mathematical Society},
  address   = {Providence, RI},
  year      = {1996},
  doi       = {10.1090/memo/0568}
}

@article{Kemp1968,
  author    = {Kemperman, J. H. B.},
  title     = {The general moment problem, a geometric approach},
  journal   = {Annals of Mathematical Statistics},
  volume    = {39},
  number    = {1},
  pages     = {93--122},
  year      = {1968},
  doi       = {10.1214/aoms/1177698508}
}

@book{KN1977,
  author    = {Kre{\u\i}n, M. G. and Nudel'man, A. A.},
  title     = {The {M}arkov moment problem and extremal problems},
  series    = {Translations of Mathematical Monographs},
  volume    = {50},
  publisher = {American Mathematical Society},
  address   = {Providence, RI},
  year      = {1977},
  note      = {Translated from the Russian (1973)}
}

@article{BarndorffNielsenPollardShephard2012,
  author    = {Barndorff-Nielsen, O. E. and Pollard, D. G. and Shephard, N.},
  title     = {Integer-valued {L}\'evy processes and low latency financial econometrics},
  journal   = {Quantitative Finance},
  volume    = {12},
  number    = {4},
  pages     = {587--605},
  year      = {2012},
  doi       = {10.1080/14697688.2012.664935}
}

@article{Freeland2010,
  author    = {Freeland, R. K.},
  title     = {True integer value time series},
  journal   = {AStA Advances in Statistical Analysis},
  volume    = {94},
  number    = {3},
  pages     = {217--229},
  year      = {2010},
  doi       = {10.1007/s10182-010-0135-0}
}

@article{Irwin1937,
  author    = {Irwin, J. O.},
  title     = {The frequency distribution of the difference between two independent
             variates following the same {P}oisson distribution},
  journal   = {Journal of the Royal Statistical Society},
  volume    = {100},
  number    = {3},
  pages     = {415--416},
  year      = {1937}
}

@article{KarlisNtzoufras2009,
  author    = {Karlis, D. and Ntzoufras, I.},
  title     = {Bayesian modelling of football outcomes: using the {S}kellam's
             distribution for the goal difference},
  journal   = {IMA Journal of Management Mathematics},
  volume    = {20},
  number    = {2},
  pages     = {133--145},
  year      = {2009}
}

@article{Maher1982,
  author    = {Maher, M. J.},
  title     = {Modelling association football scores},
  journal   = {Statistica Neerlandica},
  volume    = {36},
  number    = {3},
  pages     = {109--118},
  year      = {1982},
  doi       = {10.1111/j.1467-9574.1982.tb00782.x}
}

@article{Prekopa1990,
  author    = {Pr\'ekopa, A.},
  title     = {The discrete moment problem and linear programming},
  journal   = {Discrete Applied Mathematics},
  volume    = {27},
  number    = {3},
  pages     = {235--254},
  year      = {1990},
  doi       = {10.1016/0166-218X(90)90068-N}
}

@article{BertsimasPopescu2005,
  author    = {Bertsimas, D. and Popescu, I.},
  title     = {Optimal inequalities in probability theory: a convex optimization approach},
  journal   = {SIAM Journal on Optimization},
  volume    = {15},
  number    = {3},
  pages     = {780--804},
  year      = {2005},
  doi       = {10.1137/S1052623401399903}
}

@article{CurtoFialkow1991,
  author    = {Curto, R. E. and Fialkow, L. A.},
  title     = {Recursiveness, positivity, and truncated moment problems},
  journal   = {Houston Journal of Mathematics},
  volume    = {17},
  number    = {4},
  pages     = {603--635},
  year      = {1991}
}

@book{Grunbaum2003ConvexPolytopes,
  author    = {Gr{\"u}nbaum, B.},
  title     = {Convex polytopes},
  edition   = {Second},
  series    = {Graduate Texts in Mathematics},
  volume    = {221},
  publisher = {Springer},
  address   = {New York},
  year      = {2003},
  doi       = {10.1007/978-1-4613-0019-9},
  note      = {Prepared by V. Kaibel, V. Klee, and G. M. Ziegler}
}

@article{Hamburger1920,
  author    = {Hamburger, H.},
  title     = {{\"U}ber eine {Erweiterung} des {Stieltjes}schen {Momentenproblems}},
  journal   = {Mathematische Annalen},
  volume    = {81},
  number    = {2--4},
  pages     = {235--319},
  year      = {1920},
  doi       = {10.1007/BF01564869},
  language  = {German}
}

@article{InfusinoKunaLebowitzSpeer2017,
  author    = {Infusino, M. and Kuna, T. and Lebowitz, J. L. and Speer, E. R.},
  title     = {The truncated moment problem on {$\mathbb{N}_0$}},
  journal   = {Journal of Mathematical Analysis and Applications},
  volume    = {452},
  number    = {1},
  pages     = {443--468},
  year      = {2017},
  doi       = {10.1016/j.jmaa.2017.02.060}
}

@book{KarlinStudden1966,
  author    = {Karlin, S. and Studden, W. J.},
  title     = {Tchebycheff systems: with applications in analysis and statistics},
  series    = {Pure and Applied Mathematics},
  volume    = {15},
  publisher = {Interscience Publishers},
  address   = {New York},
  year      = {1966}
}

@book{Lang2002,
  author    = {Lang, S.},
  title     = {Algebra},
  edition   = {Third},
  series    = {Graduate Texts in Mathematics},
  volume    = {211},
  publisher = {Springer},
  address   = {New York},
  year      = {2002},
  doi       = {10.1007/978-1-4613-0041-0}
}

@article{Richter1957,
  author    = {Richter, H.},
  title     = {{Parameterfreie Absch{\"a}tzung und Realisierung von Erwartungswerten}},
  journal   = {Bl{\"a}tter der Deutschen Gesellschaft f{\"u}r Versicherungsmathematik},
  volume    = {3},
  number    = {2},
  pages     = {147--162},
  year      = {1957},
  doi       = {10.1007/BF02808864},
  language  = {German}
}

@article{Rogosinski1958,
  author    = {Rogosinski, W. W.},
  title     = {Moments of non-negative mass},
  journal   = {Proceedings of the Royal Society of London. Series A, Mathematical and Physical Sciences},
  volume    = {245},
  number    = {1240},
  pages     = {1--27},
  year      = {1958},
  doi       = {10.1098/rspa.1958.0062}
}

@book{Schmudgen2017,
  author    = {Schm{\"u}dgen, K.},
  title     = {The moment problem},
  series    = {Graduate Texts in Mathematics},
  volume    = {277},
  publisher = {Springer},
  address   = {Cham},
  year      = {2017},
  doi       = {10.1007/978-3-319-64546-9}
}

@book{ShohatTamarkin1943,
  author    = {Shohat, J. A. and Tamarkin, J. D.},
  title     = {The problem of moments},
  series    = {Mathematical Surveys},
  volume    = {1},
  publisher = {American Mathematical Society},
  address   = {Providence, RI},
  year      = {1943}
}

@article{Skellam1946,
  author    = {Skellam, J. G.},
  title     = {The frequency distribution of the difference between two {Poisson} variates belonging to different populations},
  journal   = {Journal of the Royal Statistical Society},
  volume    = {109},
  number    = {3},
  pages     = {296},
  year      = {1946},
  doi       = {10.1111/j.2397-2335.1946.tb04670.x}
}

@article{Tchakaloff1957,
  author    = {Tchakaloff, V.},
  title     = {Formules de cubatures m{\'e}caniques {\`a} coefficients non n{\'e}gatifs},
  journal   = {Bulletin des Sciences Math{\'e}matiques},
  volume    = {81},
  pages     = {123--134},
  year      = {1957},
  language  = {French}
}

\end{document}